\documentclass[twocolumn]{autart}
\usepackage{amsmath}
\usepackage{graphicx}
\usepackage{subfigure}
\usepackage{cite}
\usepackage{color,soul}
\usepackage[pdftex, colorlinks=true,citecolor=blue,linkcolor=red,anchorcolor=green]{hyperref}

\newenvironment{IEEEproof}{{\bf \em{ Proof:}}\vspace{0em}}{}
\newtheorem{theorem}{\bf Theorem}{}{}
\newtheorem{lemma}{\bf Lemma}{}{}

\newtheorem{remark}{Remark}

\newtheorem{definition}{\bf Definition}{}{}

\newcommand{\real}{\mathcal{R}}

\graphicspath{{figures/}} \DeclareGraphicsExtensions{.pd,.epsf}

\begin{document}

\begin{frontmatter}
\title{Cooperative output regulation of multi-agent { network systems with dynamic edges}\thanksref{f1}}

\thanks[f1]{{This work} is supported by a grant from the National Natural Science Foundation of China (61104149, 61374174), the Zhejiang Province Natural Science Fund (LY13F030001, LY14F030003), and the Program for New Century Excellent Talents in University (NCET-11-0459).}

\author[JX]{Ji Xiang}\ead{jxiang@zju.edu.cn},    
\author[LYJ]{Yanjun Li}\ead{liyanjun@zucc.edu.cn},  
\author[DJH]{David J. Hill}\ead{dhill@eee.hku.hk},

\address[JX]{Department of System Science and
Engineering, College of Electrical Engineering, Zhejiang University,
P. R. China.}  
\address[LYJ]{School of Information and Electrical Engineering, Zhejiang University City College, P. R. China}             
\address[DJH]{Department of Electrical and Electronic Engineering, the University of Hong Kong and the school of Electrical and Information Engineering, the University of Sydney}

\begin{abstract}
\indent This paper investigates a new class of linear multi-agent network systems, in which nodes are coupled by dynamic edges in the sense that each edge has a dynamic system attached as well. The outputs of the edge dynamic systems form the external inputs of the node dynamic systems, which are termed ``neighboring inputs'' representing the coupling actions between nodes. The outputs of the node dynamic systems are the inputs of the edge dynamic systems. Several cooperative output regulation problems are posed, including output synchronization, output cooperation and master-slave output cooperation. Output cooperation is specified as making the neighboring input, a weighted sum of edge outputs, track a predefined trajectory by cooperation of node outputs. Distributed cooperative output regulation controllers depending on local state and neighboring inputs are presented, which are designed by combining feedback passivity theories and the internal model principle. A simulation example on the cooperative current control of an electrical network illustrates the potential applications of the analytical results.
\end{abstract}

\begin{keyword}
Multi-agent, dynamic edge, output cooperation, output synchronization, feedback passivity, internal model principle, electrical network.
\end{keyword}

\end{frontmatter}

\section{Introduction}
With the popularity of intelligent devices and the fast development of communication technology, multi-agent network systems (or the closely related subject of complex dynamic networks) have attracted more and more attentions in the control literature during the past decade \cite{PogromskyTCASI2001, LiuTAC2006, ZhouTCASI2006, XieIJRNC2007, YuAU2009, TianAu2009, ZhangTAC2011, ZhaoAU2011, GuanAUT2012, LiuTAC2012, SiljakBOOK2012}, {because interactions and cooperations between units become increasingly important}. Such a network system is often described by a graph, where nodes represent the dynamic subsystems and edges the interactions between these subsystems. One significant feature of these systems is that they can achieve some collective behaviours, such as synchronization, swarming, formations and so on, with each node running a decentralized or distributed feedback controller, rather than a centralized controller.

Among these collective behaviors, consensus and synchronization are the most extensively
studied ones. The term consensus {arose to mean} that all the agents have variables of interest converge to {one common value}. Since the seminal work \cite{JadbabaieTAC2003}, where conditions were presented for consensus of undirected multi-agent systems with first-order integrators, many significant results have been reported for first-order or second-order multi-agent systems. Readers can refer to recent surveys \cite{CaoTII2013, ChenCSM2013} and earlier surveys \cite{Olfati-SaberPIEEE2007, RenICS2007} for details. The term synchronization arose to mean that all the agents have variables converge to one time-varying trajectory, a common behavior in both time and space. In fact, synchronization has a long history study in the field of physics, including phase synchronization, limit-cycle synchronization and chaos synchronization. Most results are for nonlinear identical systems, but recent results have extended to non-identical nodes with relaxed synchronization concepts \cite{ZhaoAU2011, ZhaoTCASI2011, MontenbruckAUT2015}. For nonlinear synchronization in networks of dynamic systems, refer to survey \cite{ArenasPR2008}. There are several works reported for synchronization of   linear multi-agent systems with focuses on the output feedback in recent years \cite{TunaTAC2009, SeoAU2009, ScardoviAU2009, LiTCASI2010}.
 
Whereas consensus or synchronization requires that agents are going to have an identical states, output synchronization might occur with non-identical agents and is often more interesting. Xiang et. al. studied output synchronization in networks of identical agents by using the output regulation method \cite{XiangTAC2009}. Kim et. al. studied output synchronization in networks of single-input and single-output non-identical agents \cite{KimTAC2011}. Wang et. al. presented an internal model controller for output synchronization for more general heterogeneous multi-agents systems \cite{WangTAC2010}. It is proved in \cite{WielandAU2011} that the internal model principle is a necessary and sufficient condition for non-trivial output synchronizations. Grip et. al. studied the output synchronization problem of general right-invertible linear node systems with no knowledge about their own state or output but there is a knowledge of the relative outputs \cite{GripAUT2012}. An almost output synchronization was addressed in \cite{PeymaniAUT2014}, where the output synchronization error due to the disturbances is 
optimized in terms of the $H_\infty$ norm. 

In most studies on consensus, agents do not have interactions with each other before their controllers are added. The received or measured neighboring information forms virtual edges between agents.  {Each edge} can be thought of as a simple algebraic map to get a relative error between two connected agents. However, there are many real large systems in which their subsystems are inherently coupled to each other, such as power networks, ecological systems and so on. Such a system was seldom described by a graph and studied by graph theories, partly because the network structure is not gotten enough attention. Previous studies on these  {large} systems focus on stability analysis or decentralized controllers whose purpose is to overcome the coupling influences on stability \cite{DavisonTAC1976a, MoylanTAC1978, ArakiTAC1978, MichelTAC1983}.

This paper presents a new class of inherently coupled multi-agent network systems that has both dynamic nodes and dynamic edges. Outputs of edge dynamic systems combine to form external inputs of node dynamic systems, which are termed {\em neighboring inputs}; while the outputs of node dynamic systems are the inputs of edge dynamic systems. Several cooperative output regulation problems are studied, including output synchronization, output cooperation and master-slave output cooperation. Output cooperation in this study means that the nodes have their outputs cooperate for some objective that is specified as making the neighboring inputs track some reference trajectories. The proposed controller is distributed in the sense that the feedback information contains not only the local state of the agent itself, but also the neighboring input which {contains} some indirect information of neighboring agents. There are several works about adaptively adjusting the edge weights \cite{LiTAC2015, ShafiTCNS2014}. The dynamics on the edge weights is different from those considered in this paper. In \cite{BurgerSNCS2013}, a fairly related but different work was reported. There the agents interact with each other by the controllers placed on edges, so it is the edge dynamics rather than the node dynamics to be designed to achieve output synchronization.

The development here is passivity-based. The edge dynamics is assumed to {be} strictly passive and the node dynamics to be feedback passive from the neighboring input to the nodal local output. There are several works which exploit passivity for seeking consensus or synchronization of multi-agent systems. In \cite{ChopraBK2006}, output synchronization in networks of nonlinear systems that are input-output passive was investigated. In \cite{ArcakTAC2007}, a passivity-based design is proposed for a coordination problem of second-order multi-agent systems by making the feedback channel to be passive. These results are not applicable here, because of non-Hurwitz exosystems being considered (due to the internal model principle). The closed-loop system is no longer formed by a negative feedback interconnection of two passive systems.

The main originality of this paper is two fold. One is the new model with dynamic edges and the output cooperation problem thereof. The other is the idea combining the tools of passification and internal model, and the solution based on it. The first tool together with the passivity of edge dynamics leads to a decentralized way to solve the problems and the second tool ensures a no-bias trajectory tracking.

The remainder of this paper is organized as follows: Section \ref{sec02} presents the new multi-agent network system with dynamic edges and specifies two problems of output cooperation and output synchronization. Section \ref{sec03} addresses the output synchronization problem by three subsections of internal model control, passification design and solution for output synchronization. Sections \ref{sec04} and \ref{sec05} address the output cooperation and the master-slave output cooperation problems, respectively. Section \ref{sec06} provides a numerical simulation showing applications of the developed results on cooperative current control of an electrical network. Section \ref{sec07} concludes this paper. All the proofs are placed in the Appendix. {A brief version without proofs has been presented in} \cite{XiangIFAC2014}.

\begin{figure*}[t] 
  \centering 
  \subfigure[Using dynamic edges]{ 
    \label{fig01:a} 
    \includegraphics[width=0.455\textwidth]{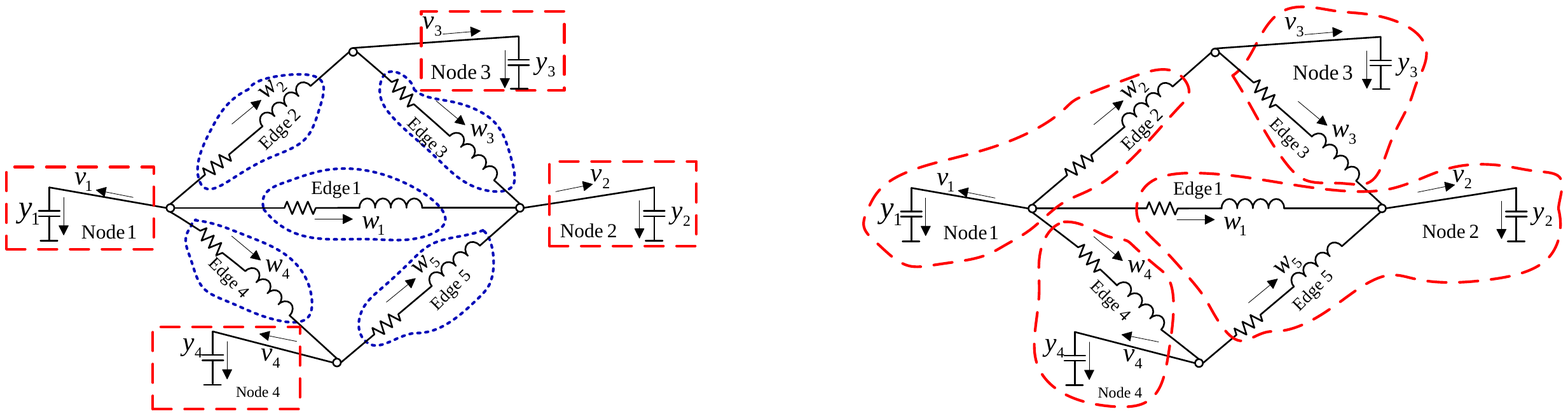}} 
  \hfill
  \subfigure[Not using dynamic edges]{ 
    \label{fig01:b} 
    \includegraphics[width=0.465\textwidth]{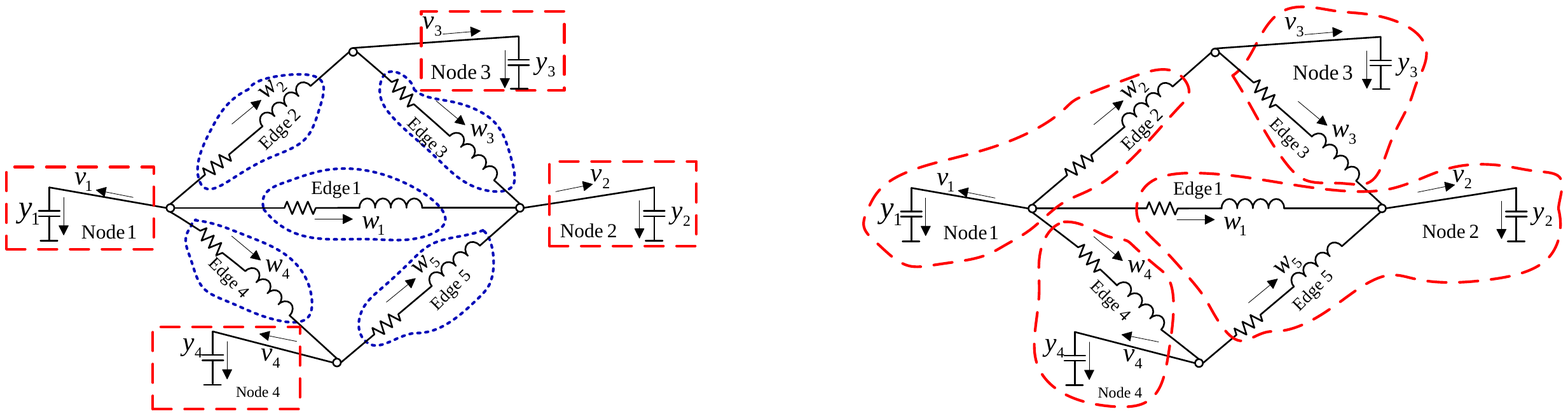}}
  \caption{A motivation example of dynamic edges} 
  \label{fig01} 
\end{figure*}

\section{Problem Formulation} \label{sec02}
Consider a multi-agent system of $N$ nodes and $M$ edges, where the node dynamics have the form of
\begin{equation}
  \label{II01}
  \left\{\begin{split}
  &\dot{x}_i = A_i x_i + B_iu_i + D_i v_i \\
  & y_i = C_i x_i
  \end{split}\right., \quad i = 1,2,\cdots, N,
\end{equation}
where $x_i\in\real^{n_i}$ is the state of node $i$, $u_i\in\real^{m_i}$ the input, $y_i\in\real^{p}$ the output and $v_i\in\real^{p}$ is the neighboring input to represent influences from other nodes. $A_i$, $B_i$, $C_i$ and $D_i$ are constant matrices with compatible dimensions.

Differently from previous formulations where $v_i$ is an algebraic function of $x_i$ and $x_j$ with node $j$ being a neighboring node of node $i$, we consider here that all the nodes are coupled by dynamic edges, that is, each edge has a dynamic system model as well,
\begin{equation} \label{II02}
\left\{
\begin{split}
  &\dot{z}_{i} = E_i z_{i} + F_i s_{i}\\
  & w_i = G_i z_i
  \end{split}\right., \quad i=1,2,\cdots,M,
\end{equation}
where $z_i\in\real^{n_{ei}}$, $s_i\in\real^{m_{ei}}$ and $w_i\in\real^{p}$ are the state, input and output of edge $i$, respectively. $E_i$, $F_i$ and $G_i$ are constant matrices with compatible dimensions. 
\begin{description}
\item [A1)] Matrices $B_i$, $D_i$ and $F_i$ are assumed to be of full column rank; $C_i$ and $G_i$ are assumed to be of full row rank.
\end{description}
The $N\times M$ incidence matrix $H$ describes the coupling relationship between the nodes, and is defined as
\begin{equation} \label{II03}
  h_{ij} = \left\{
  \begin{array}{ll}
  +1, & \mbox{node $i$ is at the positive end of edge $j$},\\
  -1, & \mbox{node $i$ is at the negative end of edge $j$},\\
  0, & \mathrm{otherwise.}
  \end{array}
  \right.
\end{equation}
The orientation of each edge only reflects that {the influence} of an edge on two nodes connected by the edge {are opposite} and can be set arbitrarily. The node neighboring input $v_i$ and the edge input $s_i$ are assumed to satisfy
\begin{equation} \label{II04}
\left\{
\begin{split}
  &v_i = -\sum_{j=1}^M h_{ij}w_j, \quad i=1,2,\cdots, N,\\
  &s_i = \sum_{j=1}^N h_{ji} y_j,\quad i=1,2,\cdots,M
  \end{split}\right..
\end{equation}
The dynamic edges model above, as well as the cooperative output regulation problem to be defined later, can be illustrated using an electrical power network that as a strongly nonlinear interconnected system, is very difficult for analysis and control. Fig. \ref{fig01} illustrates a simplified power network, where four generator nodes {have output voltage $y_i$ for the local load, and are coupled to} each other by the transmission line that is a dynamic system due to the presence of inductance. Neighboring input $v_i$ 
influencing node output $y_i$ is formed by the edge outputs $w_j$ that is in turn determined by the node output $y_i$. Therefore, both regulating the node output and regulating the neighboring input are cooperative behaviors between nodes. 

The network in Fig. \ref{fig01} can also be described by a multi-agent model without dynamic edges by incorporating the dynamic edges into the node systems. Clearly there are multiple ways to do this. Fig. \ref{fig01:b} provides an example, where edge 1 and edge 5 are placed into the node 2 and the remained edges are placed into nodes one to one. Such a model has at least two drawbacks, compared with the dynamic edge model in Fig. \ref{fig01:a}. One is that incorporation of edges losses structural information of edges and moreover {the optimal incorporation is difficult to find}. The other is that the coupling is not uniform between nodes.

Our goal is to design $u_i$ for each node $i$, $i=1,\cdots, N$, which depends on the information of  $x_i$ and $v_i$, so as to cooperatively regulate the neighboring inputs to track some trajectories given by
\begin{equation} \label{II05}
  \dot{\nu}_{i} = S_\eta \nu_i, \quad \bar{v}_i= Q_v \nu_i, \quad i=1,\cdots, N,
\end{equation}
where $\nu_i\in\real^q$ is the state of neighboring input reference system, $Q_v\in\real^{p\times q}$ is the output matrix, and $S_\eta\in\real^{q\times q}$ satisfies
\begin{description}
\item [A2)] Matrix $S_\eta$ has all eigenvalues locate on the imaginary axis with the algebraic multiplicity of $1$. 
\end{description}

\begin{remark}
Assumption A2) means that there is a symmetric positive definite matrix $P_\eta$ such that $P_\eta S_\eta + S_\eta^T P_\eta\le 0$.
\end{remark}

Formally, the following problem is proposed,

\textbf{Output Cooperation Problem:}
Given a multi-agent system consisting of dynamic systems \eqref{II01} and \eqref{II02} with the relationships \eqref{II03} and \eqref{II04}, design a distributed control law depending on local states $x_i$ and neighboring inputs $v_i$ such that the closed-loop system has each neighboring input asymptotically converge to its reference trajectory given by \eqref{II05}, namely,
\begin{equation}
v_i\rightarrow \bar{v}_i,\quad \mbox{for all } \,i=1,2,\cdots,N. 
\end{equation}
There are many real scenarios requiring the output cooperation, among which regulating neighboring input, as illustrated in Fig. \ref{fig01}, corresponds to the output current control of generators. Cooperatively pushing an elastic object along the predefined trajectory by mobile robots is another example of output cooperation, since the object has dynamics due to the elastic contact with robots. The output cooperation is also valid for the conventional multi-agent system without dynamic edges. In such a case, $w_i=G'_iF_is_i$ for some matrix $G'_i$ and the goal becomes making some weighted sum of nodal outputs track given trajectories.

One of the challenging points of the output cooperation problem is that the variable to be regulated is not a node state variable but a combination of outputs of edge dynamic systems that are driven by the node outputs. 

Notice that the node outputs are the inputs of edge systems. According to the internal model principle it is necessary for output cooperation problem that the node output has the mode of $S_\eta$. In this consideration, we propose another cooperative output regulation problem synchronizing the node outputs to a common trajectory $y_\eta$ that is given by
\begin{equation} \label{II07}
  \dot{\eta}_0=S_\eta \eta_0, \quad y_\eta = Q_\eta \eta_0,
\end{equation}
where $\eta_0\in\real^{q}$ is the exosystem state, $y_{\eta }\in\real^p$ the exosystem output, and $Q_\eta$ the output matrix. 

\textbf{Output Synchronization Problem:} Given a multi-agent system consisting of dynamic systems \eqref{II01} and \eqref{II02} with the relationships \eqref{II03} and \eqref{II04}, design a distributed control law depending on  {local states $x_i$} and neighboring inputs $v_i$ such that all the nodes have their outputs in the closed-loop system asymptotically converge to a nontrivial common trajectory $y_\eta$.

\begin{remark} Output synchronization considered here not only is of important significance by itself, but also, as seen later, plays a key step for developing results of the output cooperation problem.
It should be pointed out that $y_\eta$ defined above is a family of trajectories given that the initial condition $\eta_0(0)$ is arbitrary. This means that $y_\eta$ is not known before output synchronization; whereas $\bar{v}_i$ in the output cooperation problem are predefined and known all the time.
\end{remark}
An assumption for the edge dynamic system is made as follows,
\begin{description}
\item [A3)] The edge dynamic system  $(E_i,F_i,G_i)$ is strictly positive real in the sense defined in \cite{WenTAC1988}.
\end{description}
In the output synchronization state, $s_i=0$, subsequently $w_j=0$ since $E_i$ is Hurwitz, and then $v_i=0$. Therefore, the output synchronization problem can be regarded as a special case of output cooperation problem with $\nu_i(0)=0$. 

The network topology is assumed to satisfy
\begin{description}
\item [A4)] The network is connected, namely, the row rank of the incidence matrix $H$ is $N-1$.
\end{description}
This assumption implies that only ${\underline{\mathbf{1}}}$ is the null space base of $H^T$, i.e., $H^T\underline{\mathbf{1}}=0$, where $\underline{\mathbf{1}}$ denotes the vector with all elements being $1$. 
Define a matrix $T\in\real^{(N-1)\times N}$ satisfying $T\underline{\mathbf{1}}=0$ and $TT^T=I_{N-1}$ and an induced matrix $\bar{H}=TH$. It can be verified that $\bar{H}$ is of full row rank. Matrices $\bar{H}$ and $T$ will be often used in the rest of this paper. 

Before to end this section, let us return to the electrical network in Fig. \ref{fig01}. If $S_\eta=[0,w;-w,0]$ for some angle frequency $w$, then the {\em output cooperation} and {\em output synchronization} problems are to make the output currents track some given sinusoid and to make the output voltages not only synchronize but also be a sinusoid with frequency $w$, respectively. Although the most difficult frequency synchronization problem is avoided by setting a common $S_\eta$, the development below shows that the analysis and control is still a challenging problem even for an electrical network that is linear when the frequency is fixed and known. Assumptions for node dynamics are given in the next section along with the controller design. 

\section{Output synchronization} \label{sec03}
\subsection{Internal model controller}
Since all the node outputs should synchronize on a trajectory that is unknown but determined by $(S_\eta,Q_\eta)$, see \eqref{II07}, a natural starting point is that each node makes its output track a reference trajectory, which is independent of others and produced by 
\begin{equation} \label{III01}
  \dot{\eta}_i = S_\eta\eta_i , \quad y_{\eta i}=Q_\eta\eta_i, \quad i=1,2,\cdots,N,
\end{equation}
where $\eta_i\in\real^{q}$ is the exosystem state, and $y_{\eta i}\in\real^p$ the exosystem output. The tracking error is defined by
\begin{equation} \label{III02}
  e_i=y_i-y_{\eta i}=C_ix_i - Q_\eta\eta_i.
\end{equation}
An internal model controller to make $e_i\rightarrow 0$ has the form of
\begin{equation} \label{III03}
  \left\{
  \begin{split}
  &\dot{\zeta}_i  = G_{1 i} \zeta_i + G_{2 i}(y_i-Q_\eta\eta_i)\\
  & u_i = K_{xi}x_i + K_{\zeta i}\zeta_i
  \end{split}
  \right.
\end{equation}
where $\zeta_i\in\real^{c_i}$ is the controller state, matrix pair $(G_{i1}, G_{i2})$ incorporates a $p$-copy internal model of matrix $S_\eta$, and $K_{xi}$ and $K_{\eta i}$ are feedback gains to be designed in the next subsection. 

The $p$-copy internal model is a crucial skill to address robust output regulation problem with $p$ dimensional outputs , which is recalled as follows for readability \cite{HuangBOOK2004}.

\begin{definition}[$p$-copy internal model]
Given a matrix $S_\eta$, a pair $(G_1,G_2)$ is said to incorporate a $p$-copy internal model of $S_\eta$ if the pair $(G_1,G_2)$ admits the following form
\begin{equation}
G_1 = T_p \begin{bmatrix} S_1 & S_2 \\ 0 & G_{p1} \end{bmatrix}T^{-1}_p, \quad 
G_2 = T_p \begin{bmatrix} S_3 \\ G_{p2} \end{bmatrix},
\end{equation}
where $S_1, S_2, S_3$ are arbitrary matrices with compatible dimensions, $T_p$ is any non-singular matrix with same dimension as $G_1$ and $(G_{p1}, G_{p2})$ is described as follows
\begin{equation}
G_{p1}=\mathrm{diag} \underbrace{[\alpha_1,\cdots, \alpha_p]}_{\mathrm{p-tuple}},
G_{p1}=\mathrm{diag} \underbrace{[\beta_1,\cdots, \beta_p]}_{\mathrm{p-tuple}},
\end{equation} 
where for $i=1,\cdots,p$, $\alpha_i$ is a constant square matrix of dimension $d_i$ for some integer $d_i$ and $\beta_i$ is a constant column vector of dimension $d_i$ such that
\begin{description}
\item [(i)] $\alpha_i$ and $\beta_i$ are controllable.
\item [(ii)] The minimal polynomial of $S_\eta$ divides the characteristic polynomial of $\alpha_i$.
\end{description}
\end{definition}
Here for $S_\eta$ satisfying Assumption A2), the minimal polynomial is the same as the characteristic polynomial and $d_i=q$. Let $a_q(\lambda)=\lambda^q+a_1\lambda^{q-1}+\cdots+a_{q-1}\lambda+a_q$ be the minimal polynomial of $S_\eta$. A general way to show how the internal model forces the tracking error to be zero is to select $(\alpha_i,\beta_i)$ in the controllable canonical form described by
\begin{equation}
\alpha_i=
\begin{bmatrix}
0 & 1 & \cdots & 0 \\ \vdots & \vdots & \ddots & \vdots \\
0 & 0 & \cdots & 1 \\ -a_q & -a_{q-1} & \cdots & -a_1
\end{bmatrix}, \quad \beta_i=\begin{bmatrix} 0 \\ \vdots \\ 0 \\ 1 \end{bmatrix}.
\end{equation}
Taking $p=1$ for example, the first equation in \eqref{III03} becomes 
\begin{equation} \label{rev01}
\dot{\zeta}_i=\alpha_i\zeta+\beta(y_i-Q_\eta\eta_i).
\end{equation}
When the closed-loop system is exponentially stable, the system steady state is driven by the exosystem, namely, there are matrices $\Theta$ and $\Pi$ such that $\zeta_i\rightarrow \Theta \eta_i$ and $x_i\rightarrow \Pi \eta_i$. Substituting $\zeta_i=\Theta\eta_i$ into \eqref{rev01} yields $\theta_2=\theta_1 S_\eta, \theta_3=\theta_1 S_\eta^2,\cdots, \theta_q=\theta_1S^{q-1}_\eta$ and $\theta_qS_\eta+\sum_{i=1}^q{a_i}\theta_{q-i+1}=y_i-Q_\eta\eta_i$, where $\theta_i$ denotes the $i$th row of $\Theta$. The left side of the last equality is nothing but $\theta_1(S_\eta^d+a_1S_\eta^{d-1}+\cdots+a_d)=0$ due to the Cayley-Hamilton theorem. This means that $y_i\rightarrow Q_\eta \eta_i$ and $C_i\Pi=Q_\eta$. In summary, the Lemma 1.27 in \cite{HuangBOOK2004} is simplified as follows,
\begin{lemma}\label{le01a}
The zero tracking error is guaranteed by the $p$-copy internal model controller \eqref{III03} if the resulted closed-loop system is exponentially stable.
\end{lemma}
 
\subsection{Passification design}
With controller \eqref{III03}, the dynamics of node system \eqref{II01} becomes
\begin{equation}
  \label{III04}
\dot{\hat{x}}_i = \hat{A}_i\hat{x}_i + \hat{D}_i v_i + \hat{D}_{\eta i}\eta_i, \quad y_i =\hat{C}_i\hat{x}_i
\end{equation}
where $\hat{x}_i = [x_i^T, \zeta_i^T]^T$,
\begin{gather*}
\hat{A}_i = \begin{bmatrix}
  A_i + B_iK_{x i} & B_i K_{\zeta i} \\
  G_{2i}C_i & G_{1 i}
\end{bmatrix}, \quad \hat{D}_i = \begin{bmatrix}
  D_i \\ 0
\end{bmatrix}, \\ \hat{D}_{\eta i} = \begin{bmatrix}
  0 \\ -G_{2 i }Q_\eta
\end{bmatrix}, \quad \hat{C}_i =\begin{bmatrix}
  C_i & 0
\end{bmatrix}.
\end{gather*}
In a single-node system, it is enough to make $\hat{A}_i$ be Hurwitz for output regulation. Here for cooperation among nodes, an extra passification is required for the design of feedback gains $K_{xi}, K_{\zeta i}$. 
\begin{description}
\item [A5)] The feedback gains $K_{xi}$ and $K_{\zeta i}$ are designed in such a way that  $\hat{A}_i$ is Hurwitz and the closed-loop system $(\hat{A}_i,\hat{D}_i, \hat{C}_i)$ is passive.
\end{description}
Making $\hat{A}_i$ Hurwitz requires that node $i$ satisfies the following assumptions,
\begin{description}
\item [A51)] Matrix pair $(A_i,B_i)$ is stabilizable.
\item [A52)] For any eigenvalue $\lambda$ of $S_\eta$,\\[-0.5em]
\[
\mbox{rank}\left(\begin{bmatrix}
  A_i-\lambda I & B_i \\ C_i & 0
\end{bmatrix}\right)=n_i+p.
\]
\end{description}
There has been plentiful works on the passification design (see \cite{BarkanaTAC2004} \cite{AndrievskiiARC2006} for surveys of this area) to make a system be a "direct" passivation with respect to input, but seldom works are reported for rendering a system an ``indirect'' passivation with respect to the disturbances input. Arcak and Kokotovi\'{c} \cite{ArcakSCL2001} presented the feasible condition for rendering a single-input-single-output (SISO) system strictly positive real with respect to the disturbance input. The feasibility condition for the indirect passification design required in Assumption A5) is still open problem. However, if the node system is restricted to the direct passification case, i.e, $B_i=D_i$, one feasibility condition can be stated as follows, 
\begin{description}
\item [A53)] $C_iB_i$ is symmetric positive definite (SPD). 
\item [A54)] The polynomial $\varphi_0(\lambda)=\det\left[\begin{smallmatrix}
  A_i-\lambda I & B_i \\ C_i & 0
\end{smallmatrix}\right]$ is Hurwitz, i.e, minimum phase condition, 
\end{description}

The combination of the two assumptions is called the hyper minimum phase condition \cite{AndrievskyCDC1996}. As illustrated below, this condition is also valid for the passification by the internal mode controller \eqref{III03} with assumption A2). 
\begin{theorem} \label{th00}
Given a node system \eqref{II01} with $B_i=D_i$ and with assumptions A2), A53) and A54), there exist feedback gains $K_{xi}$ and $K_{\zeta i}$ such that the closed-loop system \eqref{III04} satisfies Assumption A5).	
\end{theorem}
 
Assumption A53) restricts that the system is of the uniform relative degree one.  Assumption A54) together with A2) implies A52). Since our main focus is not on the passification design, in this paper we use Assumption A5) instead of A51)$\sim$A54) to leave the possibility for more general node systems that could be of indirect passification.  

\subsection{Solution for output synchronization}
Since $\hat{A}_i$ and $S_\eta$ have no common eigenvalue, the following Sylvester equation has a unique solution $\Pi_i$,
\begin{equation} \label{III05}
  \Pi_iS_\eta = \hat{A}_i\Pi_i + \hat{D}_{\eta i},
\end{equation}
and, since matrix pair $(G_{i1}, G_{i2})$ incorporates a $p$-copy internal model of matrix $S_\eta$, the solution $\Pi_i$ further satisfies
\begin{equation} \label{III06}
  \hat{C}_i\Pi_i = Q_\eta,
\end{equation}
which implies that $e_i\rightarrow 0$ under controller \eqref{III03} if no coupling between nodes, i.e., $v_i=0$. In the presence of coupling, making $e_i\rightarrow 0$ is a decentralized servomechanism problem \cite{DavisonTAC1976a}, because $\eta_i$ is independent of $\hat{x}_i$ and $\zeta_i$. For this goal, the following result is given,

\begin{theorem} \label{th01}
Given a multi-agent system consisting of \eqref{II01}$\sim$\eqref{II04} with assumption A1) and exosystem \eqref{III01}, if assumptions A2), A3) and A5) hold, then $e_i$ will exponentially converge to zero for all $i=1,\cdots,N$, under controller \eqref{III03}.
\end{theorem}

The above theorem implies that when the exosystem $\eta_i$, whose output the node $i$ will track, have the same dynamic model for all the nodes, then the decentralized internal model controller can realize the output tracking for the networks coupled by dynamic edges if some passivity properties are satisfied.

\begin{figure}
  \centering
  \includegraphics[width=0.49\textwidth]{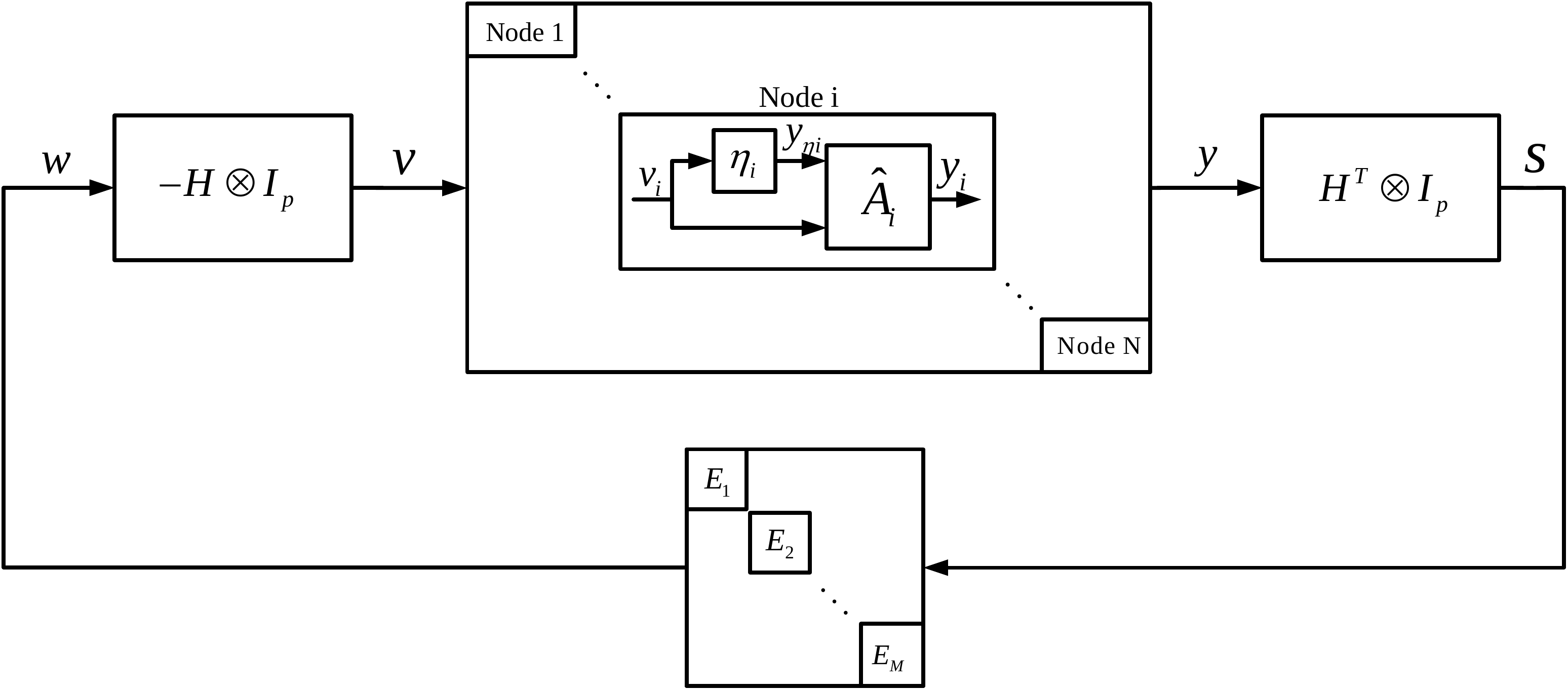}\\
  \caption{A block diagram representation for the interconnection of multi-agent systems \eqref{II01}$\sim$\eqref{II04} with controller \eqref{III07}.}\label{fig02}
\end{figure}

Now we consider the output synchronization problem. If these $\eta_i$ are synchronous to each other, then $y_i-y_j\rightarrow 0$ can be {further} obtained by controller \eqref{III03}. But in general, $\eta_i\neq \eta_j$ due to different initial conditions\footnote{In real, the same initial state can not be achieved because of among others the different staring time between nodes.}. Meanwhile, the synchronization errors, either $y_i-y_j$ or $\eta_i-\eta_j$, are not directly available for synchronization seeking of $\eta_i$. Here neighboring input $v_i$ is the only available information that indirectly reflects the synchronization error. Our idea is to adjust the exosystem dynamics \eqref{III01} by  using $v_i$ in order to synchronize $\eta_i$. With this alteration, the following dynamic controller is presented,
\begin{equation} \label{III07}
\left\{
\begin{split}
  &\dot{\eta}_i=S_\eta\eta_i + \epsilon B_{\eta} v_i \\
  &\dot{\zeta}_i  = G_{1 i} \zeta_i + G_{2 i}(y_i-Q_\eta\eta_i)\\
  & u_i = K_{x i}x_i + K_{\zeta i}\zeta_i
\end{split}\quad ,\right.
\end{equation}
where $B_{\eta}\in\real^{q\times p}$ is the input matrix of the exosystem, $\epsilon$ is a positive scalar expressing the strength of adjusting the dynamics of output reference, and $(G_{i1}$, $G_{i2})$, $K_{xi}$ and $K_{\zeta i}$  are the same as those in \eqref{III03}. The block diagram of the resulting closed-loop system is shown in Fig. \ref{fig02}.

\begin{theorem} \label{th02}
Given a multi-agent system consisting of \eqref{II01}$\sim$\eqref{II04} with assumption A1). If Assumptions A2)$\sim$A5) hold and $B_\eta$ is designed such that $(S_\eta, B_\eta, Q_\eta)$ is passive, namely, $B_\eta=P_\eta^{-1}Q_\eta^T$, then there is a scalar $\epsilon^*>0$ such that for all $0<\epsilon<\epsilon^*$, controller \eqref{III07} will solve the output synchronization problem. And moreover, the steady output trajectory is given by
\begin{equation} \label{III08}
  y_i(t) \rightarrow Q_\eta e^{S_\eta t}\frac{1}{N}\sum_{i=1}^N\eta_i(0), \quad i=1,2,\cdots,N.
\end{equation}
\end{theorem}

As shown in Fig. \ref{fig02}, the signal channel from reference output $y_{\eta i}$ to $s_j$ is not direct, but goes indirectly through output $y_i$ of dynamic system $\hat{x}_i$. Critically, it is {no longer passive} from $y_{\eta i}$ to $y_i$ due to $\hat{C}_i\hat{D}_{\eta i} =0$. This in turn results in the adjusting strength of an upper bound $\epsilon^*$ since the closed-loop system is not a feedback interconnection of two passive systems.

\section{Output Cooperation} \label{sec04}
This section addresses the output cooperation problem making the influences between nodes, the neighboring inputs $v_i$, track some predefined trajectories $\bar{v}_i$. Similar with the output synchronization problem, the trick is to adjust output reference $\eta_i$ according to tracking error $v_i-\bar{v}_i$. 

Suppose that the node output is just its reference, $y_i=y_{\eta_i}$ a perfect tracking, then the output reference system \eqref{III01} can be regarded as a dynamic controller for the edge dynamic systems in the sense that $y_{\eta_i}$ is the control input to regulate $v_i\rightarrow \bar{v}_i$, replacing $y_i$ with $y_{\eta_i}$. In this consideration, the dynamics of $\eta_i$ should be altered to incorporate a $p$-copy model of $S_\eta$. In summary the controller is given by,
\begin{equation} \label{IV02}
\left\{
\begin{split}
  &\dot{\bar{\eta}}_i=G_{S}\bar{\eta}_i + \epsilon G_B (v_i-Q_v\nu_i) \\
  &\dot{\zeta}_i  = G_{1 i} \zeta_i + G_{2 i}(y_i-G_Q\bar{\eta}_i)\\
  & u_i = K_{x i}x_i + K_{\zeta i}\zeta_i
\end{split},\right.
\end{equation}
where $\bar{\eta}_i\in\real^{pq}$. Controllable matrix pair $(G_S,G_B)$ incorporates a $p$-copy model of $S_\eta$. A possible choice for $G_S$, $G_B$ and $G_Q$ is
\begin{gather} 
  G_S= I_p\otimes S_\eta, \nonumber\\
  \quad G_B = \begin{bmatrix}
    B^1_{\eta} & \cdots & 0 \\
    \vdots & \ddots & \vdots \\
    0 & \cdots & B_\eta^p
  \end{bmatrix}, \quad G_Q = \begin{bmatrix}
    Q_\eta^1 & \cdots & 0 \\
    \vdots & \ddots & \vdots \\
    0 & \cdots & Q_\eta^p
  \end{bmatrix} \label{IV03}
\end{gather}
where $B_\eta^i$ denotes the $i$th column of $B_\eta$ and $Q_\eta^i$ the $i$th row of $Q_\eta$. It can be verified that
\begin{equation} \label{IV04}
  (I_p\otimes P_\eta)G_S+G_S^T(I_p\otimes P_\eta)\le 0, \,\, (I_p\otimes P_\eta)G_B = G_Q^T.
\end{equation}

\medskip
\begin{remark}
With $Q_v \nu_i \equiv 0$, controller \eqref{IV02} reduces to an output synchronization controller, especially the same as \eqref{III07} when $p=1$. In this consideration, the output synchronization problem is a special one of the output cooperation problem. Noting that an edge can be regarded as a filter of the output error between nodes that the edge connects, the regulation of neighboring inputs is regulating the weighted sum of the filtered output synchronization errors.
\end{remark}

The analysis developed below consists of three parts, the separate node dynamics, the separate edge dynamics, and the whole coupled dynamics.

Firstly, we turn back to node dynamics \eqref{III04} with $v_i$ and $\eta_i$ being replaced by $\bar{v}_i$ and $\bar{\eta}_i$, respectively,
\begin{equation} \label{IV05}
  \dot{\hat{x}}_i = \hat{A}_i \hat{x}_i  + \hat{D}_{\bar{\eta}i}\bar{\eta}_i+ \hat{D}_i Q_v \nu_i,
\end{equation}
where $\hat{D}_{\bar{\eta}_i}=\begin{bmatrix}  0 \\ -G_{2i}G_Q \end{bmatrix}$.  {The dynamics of $\hat{x}_i$ in system} \eqref{IV05}  {is driven by $\bar{\eta}_i$ and $\nu_i$.}
Since $\hat{A}_i$ is Hurwitz and $S_\eta$ has no stable eigenvalues, there is a unique solution $\bar{\Pi}_i=\begin{bmatrix}
  \bar{\Pi}_{1i}, \bar{\Pi}_{2i}
\end{bmatrix}$ with $\bar{\Pi}_{1i}\in\real^{(n_i+c_i)\times pq}, \bar{\Pi}_{2i} \in\real^{(n_i+c_i)\times q}$ satisfying
\begin{equation} \label{IV06}
  \bar{\Pi}_i\begin{bmatrix}
    G_S & 0 \\ 0 & S_\eta
  \end{bmatrix}= \hat{A}_i \bar{\Pi}_i + \begin{bmatrix}
    \hat{D}_{\bar{\eta} i} & \hat{D}_i Q_v
  \end{bmatrix}, \quad
\end{equation}
and
\begin{equation} \label{IV07}
  \hat{C}_i\bar{\Pi}_{1i} = G_Q, \quad \hat{C}_i\bar{\Pi}_{2 i}=0.
\end{equation}

Secondly, we consider the edge dynamics \eqref{II02} with a perfect output tracking of the node system; replacing $y_i$ by $G_Q\bar{\eta}_i$, together with the first equation in \eqref{IV02} and reference system \eqref{II05}, obtains the following dynamic systems
\begin{equation} \label{IV08}
\left\{
\begin{split}
  & \dot{z}_j = E_jz_j + F_j\sum_{i=1}^N h_{ij}G_Q\bar{\eta}_i,\\ 
  & \dot{\bar{\eta}}_i = G_S\bar{\eta}_i + \epsilon G_B(-\sum_{j=1}^M h_{ij}G_jz_j-Q_v\nu_i), \\
  & \dot{\nu}_i = S_\eta \nu_i, 
\end{split}\right..
\end{equation}
where $j=1,2,\cdots,M$ for the first equation and $i=1,2,\cdots,N$ for the last two equations.

Introduce a coordinates transformation $\tilde{\eta}_i=(T_i\otimes I_{pq})\bar{\eta}$, $\tilde{\nu}_i = (T_i\otimes I_{pq})\nu$, $i=1,\cdots,N-1$,  $\tilde{\eta}_N=\sum_{i=1}^N\bar{\eta}_i$, and $\tilde{\nu}_N=\sum_{i=1}^N\nu_i$, where $\bar{\eta}=[\bar{\eta}_1^T,\cdots, \bar{\eta}_N^T]^T$, $\nu=[\nu_1^T,\cdots, \nu_N^T]^T$, and $T_i$ is the $i$th row of $T$, defined after assumption A4).
Define
\begin{gather*}
z = \begin{bmatrix}
  z_1 \\\vdots \\ z_M
\end{bmatrix}, \,\, \tilde{\eta} = \begin{bmatrix}
  \tilde{\eta}_1 \\ \vdots \\ \tilde{\eta}_{N-1}
\end{bmatrix}, \,\, Z = \begin{bmatrix}
  z \\ \tilde{\eta}
\end{bmatrix},\,\,
\tilde{\nu}=\begin{bmatrix}
  \tilde{\nu}_1\\ \vdots \\ \tilde{\nu}_{N-1}
\end{bmatrix}.
\end{gather*}

With them, the system \eqref{IV08} can be transformed into
\begin{subequations}\label{IV09}
\begin{align} 
  \label{IV09a}  &\dot{Z}=\mathcal{A}_\nu Z + \mathcal{B}_\nu \tilde{\nu}, \quad \dot{\tilde{\nu}}=(I_{N-1}\otimes S_\eta)\tilde{\nu} \\ 
\label{IV09b} & \dot{\tilde{\eta}}_N = G_S\tilde{\eta}_N - \epsilon G_BQ_v \tilde{\nu}_N, \quad \dot{\tilde{\nu}}_N = S_\eta \tilde{\nu}_N
\end{align}
\end{subequations}
where
\begin{gather*}
\mathcal{A}_\nu = \begin{bmatrix}
\overline{EM} & \overline{\bar{H}^TF}(I_{N-1}\otimes G_Q) \\
-\epsilon(I_{N-1}\otimes G_B)\overline{\bar{H}G} & I_{N-1}\otimes G_S
\end{bmatrix}, \\ \mathcal{B}_\nu = \begin{bmatrix}
  0 \\ -\epsilon (I_{N-1}\otimes G_BQ_v)
\end{bmatrix},
\end{gather*}
with $\overline{EM}$ being defined in \eqref{A05}, $\overline{\bar{H}^TF}$ and $\overline{\bar{H}G}$ defined in \eqref{A15}. Here the fact that $\left[\begin{smallmatrix} T \\ \underline{\mathbf{1}}^T\end{smallmatrix}\right]\left[
\begin{smallmatrix} T^T & \underline{\mathbf{1}}/N \end{smallmatrix}
\right]=I_N$ is used.

\begin{theorem} \label{th03}
Given system \eqref{IV08} with assumptions A1)$\sim$A4) and with $\epsilon>0$, then there is a matrix $\tilde{\Pi}=[\tilde{\Pi}_z^T, \tilde{\Pi}_\eta^T]^T$ satisfying
\begin{equation} \label{IV10}
\tilde{\Pi}(I_{N-1}\otimes S_\eta)=\mathcal{A}_\nu \tilde{\Pi} + \mathcal{B}_\nu, \,\, -\overline{\bar{H}G}\tilde{\Pi}_z=I_{N-1}\otimes Q_v
\end{equation}
such that $z\rightarrow \tilde{\Pi}_z\tilde{\nu}$, $\tilde{\eta}\rightarrow \tilde{\Pi}_\eta \tilde{\nu}$ and
\begin{equation} \label{IV11}
 v_i-Q_v\nu_i\rightarrow Q_v\nu_0(t), 
\end{equation} 
where $\nu_0(t)$ is the solution of the following dynamic system,
\begin{equation*} 
  \dot{\nu_0}=S_\eta\nu_0, \quad \nu_0(0)=-\frac{1}{N}\sum_{i=1}^N \nu_i(0).
\end{equation*}
\end{theorem}

It can be seen that system \eqref{IV08} contains two separated parts, \eqref{IV09a} and \eqref{IV09b}. The former can realize output tracking with the internal model principle by itself but cannot achieve the neighboring input tracking in general even if each node has a perfect tracking performance that $y_i=G_Q\bar{\eta}$. For system \eqref{IV09b} a nonzero $\tilde{\nu}_N$ ($\tilde{\nu}_N=N\nu_0$) will leads to an unbounded $\tilde{\eta}_N$, which means that at least some $\bar{\eta}_i$ are unbounded as well. This causes that $y_i$ will converge to infinity although the edge dynamics has bounded states and bounded inputs (weighted sums of $y_i$).

To avoid these problems, the following {\em zero sum} condition is made,

\begin{description}
  \item [A6)] The trajectories of neighboring input reference systems are in the manifold $\sum_{i=1}^N \nu_i=0$.
\end{description}
\begin{remark}
Under the zero sum condition, all the neighboring input $v_i$ will asymptotically track their reference $Q_v\nu_i$, since $\nu_0(0)=0$.
\end{remark}
At last, we are ready to present our main result for output cooperation,
\begin{theorem} \label{th04}
Given a multi-agent system consisting of \eqref{II01}$\sim$\eqref{II04}. If Assumptions A1)$\sim$A6) hold, then there is a scalar $\epsilon^*>0$ such that for all $0<\epsilon<\epsilon^*$, controller \eqref{IV02} with \eqref{IV03} and \eqref{IV04} will solve the output cooperation problem. And moreover, the sum of node output $y_s=\sum_{i=1}^Ny_i$ satisfies
\begin{equation} \label{IV12}
  y_s(t)\rightarrow G_Q e^{G_S t}\sum_{i=1}^N\bar{\eta}_i(0).
\end{equation}
\end{theorem}
It can be seen that the different initial conditions of output reference systems cause different steady outputs $y_i$, although the steady neighboring inputs $v_i$ are the same for any initial condition. This in turn means that for a given output cooperation target $\bar{v}_i$, there are infinite feasible solutions for node outputs.

\section{Master-slave output cooperation} \label{sec05}
In order for a well-posed output cooperation problem, the zero sum condition of assumption A6) is required. Even a small violation of this condition will cause some reference outputs to diverge to infinity. To solve this problem, one possible way is to keep some nodes free, running without the requirements on their neighboring inputs. Such nodes are called master nodes, in the sense that they do not receive any command so that their output reference is not changed. The controller for a master node is given by \eqref{III03}. The remaining nodes are called slave nodes whose controller is designed as \eqref{IV02}. The output reference of slave node is changed to fit in with command $\bar{v}_i$.

Without loss of generality, we assume that the first $l$ nodes are slave nodes and the remaining $N-l$ nodes are master nodes, and present the following result,
\begin{theorem} \label{th05}
Given a multi-agent system consisting of \eqref{II01}$\sim$\eqref{II04}. If the nodes have controllers described by
\begin{subequations} \label{V01}
  \begin{align}
&\left\{
\begin{array}{l}
  \dot{\bar{\eta}}_i=G_{S}\bar{\eta}_i + \epsilon G_B (v_i-Q_v\nu_i) \\
  \dot{\zeta}_i  = G_{1 i} \zeta_i + G_{2 i}(y_i-G_Q\bar{\eta}_i)\\
   u_i = K_{x i}x_i + K_{\zeta i}\zeta_i
\end{array}\right., &i=1,\cdots, l, \label{V01a}
\\
&\left\{
\begin{array}{l}
  \dot{\zeta}_i  = G_{1 i} \zeta_i + G_{2 i}(y_i-Q_\eta\eta_i)\\
   u_i = K_{x i}x_i + K_{\zeta i}\zeta_i
\end{array}\right., &i=l+1,\cdots,N, \label{V01b}
\end{align}
\end{subequations}
where $0\le l\le N-1$, and both $\nu_i$ and $\eta_i$ are the state of exosystems satisfying
\begin{equation}
\left\{
\begin{split}
  &\dot{\nu}_i=S_\eta {\nu}_i, \quad i=1,2,\cdots,l\\
  & \dot{\eta}_j=S_\eta \eta_j,\quad j=l+1,\cdots,N
\end{split}\right..
\end{equation}
If assumptions A1)$\sim$A5) are satisfied, then there is a positive scalar $\epsilon^*$ such that for all $0<\epsilon<\epsilon^*$,
\begin{equation}
  v_i\rightarrow Q_v\nu_i, \quad y_j\rightarrow Q_\eta\eta_j
\end{equation}
for all $i=1,\cdots,l$ and $j=l+1,\cdots,N$.
\end{theorem}
It can be seen that the master node plays an output tracking role to make its output $y_i$ track a given reference $Q_\eta\eta_i$; while the slave node plays an output cooperation role to make its neighboring input track given reference $Q_v\nu_i$. Theorem \ref{th05} implies that if the two kinds of nodes simultaneously exist in the network, they can definitely realize their targets without the so-called zero sum condition. Recalling Theorem \ref{th01}, the number of slave nodes can reach $l=0$.
\begin{figure}[t]
  \centering
  \includegraphics[width=0.47\textwidth]{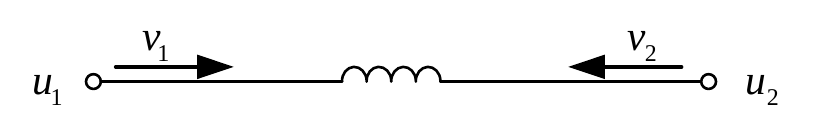}\\
  \caption{A wire connecting two voltage sources.}\label{fig03}
\end{figure}
The above master-slave configuration is illustrated by the control of microgrids working in the islanded mode \cite{GreenEPSR2007}. {One inverter, as a master generator, works in the voltage control mode to provide the fundamental voltage and other inverters, as slave generators, work in the current control mode to inject the desired powers into the microgrids. 

Here the feasibility of multiple master nodes can be explained by the example in Fig. \ref{fig03}. Two nodes can have arbitrary output voltages $u_1$ and $u_2$ but the output currents (neighboring inputs) $v_1$ and $v_2$ must satisfy the zero sum condition $v_1+v_2=0$. Therefore the case that all the nodes are master nodes, i.e., $l=0$, is feasible but the case that all the nodes are slave nodes, i.e., $l=N$, is not allowed.}

\section{Example on electrical network} \label{sec06}
\begin{figure}[!h]
  \centering
  \includegraphics[width=0.49\textwidth]{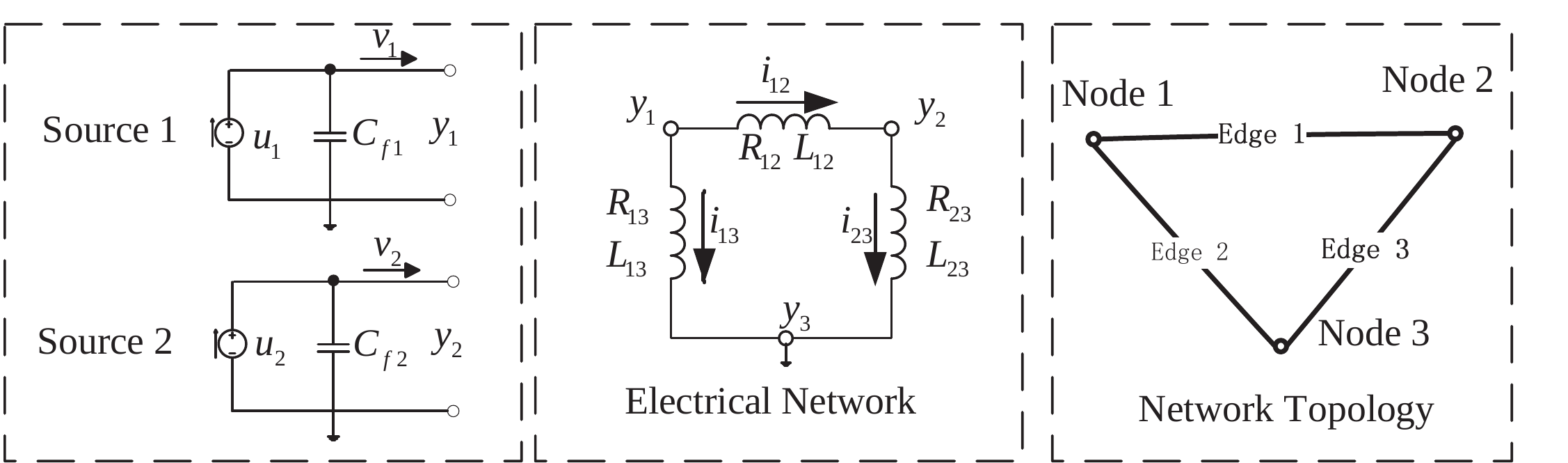}\\
  \caption{The electrical network of three nodes}\label{fig04}
\end{figure}

A simple electrical network consisting of two sources, two loads and one transmission line is selected as an application example to illustrate the analytic results. We consider the sources with their current being the control input. Fig. \ref{fig04} shows the electrical network. Node $1$ and $2$ are sources, and node $3$ is the ground. A transmission line connects the outputs of two sources, which means that both sources jointly provide currents for loads connected to node $1$ and $2$. In the model, every edge contains  a non-ideal inductor, i.e., modeled by a resistance and an inductance.

The incidence matrix is
\begin{equation}
  H=\begin{bmatrix}
    1 & 1 & 0 \\
    -1 & 0 & 1 \\
    0 & -1 & -1
  \end{bmatrix}.
\end{equation}
For the edges, the following dynamic functions can be established,
\begin{subequations}
	\begin{align}
		&L_{12}\dot{i}_{12} = -R_{12}i_{12} + y_1-y_2,\\
	&	L_{13}\dot{i}_{13}= - R_{13}i_{13}+ y_1-y_3, \\
   & L_{23}\dot{i}_{23}= - R_{23}i_{23}+ y_2-y_3.
	\end{align}
\end{subequations}
And for the sources,
\begin{subequations}
\begin{align}
	&	C_{f1}\dot{u}_{f1}=u_1-v_1,\quad y_1 = u_{f1},\\
	&	C_{f2}\dot{u}_{f2}=u_2 - v_2, \quad y_2 = u_{f2}.
\end{align}
\end{subequations}
Suppose the  desired output is a $50$hz sinewave, that is, $S_\eta =\left[\begin{smallmatrix}
  0 & -w\\
  w & 0
\end{smallmatrix}\right]$, with $w=100\pi$. Here we consider a master-slave output cooperation problem as follows. Node $3$ is the ground that can not be controlled. So it can be taken as a master node with a perfect voltage tracking performance, satisfying $y_3\equiv \eta_3$ with
\begin{equation}
  \dot{\eta}_3=S_\eta \eta_3,\quad \eta_3(0)=0.
\end{equation}
Both node $1$ and node $2$ are slave nodes, being required to make their neighboring inputs (output currents) track desired currents, $\bar{v}_1=10\sin(wt+\pi/6)$ and $\bar{v}_2=10\cos(wt)$, respectively. The desired currents are produced by the following dynamic systems,
\begin{equation}
  \dot{\nu}_i=S_\eta\nu_i,\quad  \bar{\nu_i}=Q_v\nu_i, \,\, \mathrm{with} \,\, Q_v =[1,\,\,0],
\end{equation}
where the initial conditions are
\begin{equation}  
\nu_1(0)=\begin{bmatrix}
    5 \\ -5\sqrt{3}
  \end{bmatrix}, \,\, \nu_2(0)=\begin{bmatrix}
    10 \\ 0
  \end{bmatrix}.
\end{equation}

Such a configuration of output cooperation corresponds to the scenario where two sources work in the current control mode.

Take $z_1=i_{12}$, $z_2=i_{13}$, $z_3=i_{23}$ and $x_i=u_{fi}$, $i=1,2$. Notice that
\begin{equation}
  \begin{split}
    &v_1 = i_{12}+i_{13}, \\
  & v_2 =  -i_{12}+i_{23}.
  \end{split}
\end{equation}
The electrical network is just a multi-agent system consisting of \eqref{II01}$\sim$\eqref{II04}, with assumption A1) and A4) being satisfied. The physical parameters of the  network are listed in
Table \ref{tab01}. Assumption A2) can be also verified with $P_\eta=I_2$.

\begin{table}[!hb]
\centering
  \caption{Parameters of the electrical network} 
 \begin{tabular}
    {c|c|c|c}
     \hline
    $R_{12}$ &  $L_{12}$ & $R_{13}$& $L_{13}$ \\
    \hline
   $0.05\Omega$ & $0.01$mH & $9\Omega$ &  $1$mH \\
     \hline
     $R_{23}$ & $L_{23}$ & $C_{f1}$ & $C_{f2}$  \\
     \hline
     $8\Omega$ & $5$mH & $50\mathrm{\mu F}$ & $30\mathrm{\mu F}$ \\
     \hline
  \end{tabular}
 \label{tab01}
\end{table}

According to Theorem \ref{th05}, the following controllers are designed for two source nodes,
\begin{equation}
\left\{
\begin{split}
    &\dot{\bar{\eta}}_i=G_{S}\bar{\eta}_i + \epsilon G_B (v_i-\bar{\nu}_i) \\
  &\dot{\zeta}_i  = G_{1 i} \zeta_i + G_{2 i}(y_i-G_Q\bar{\eta}_i)\\
  & u_i = K_{x i}x_i + K_{\zeta i}\zeta_i
\end{split}\right., \quad i=1,2,
\end{equation}
where
\begin{gather*}
G_S = S_\eta, \quad G_Q =\begin{bmatrix}
  0 & 1
\end{bmatrix}, \quad   G_B=  P_\eta G_Q^T=\begin{bmatrix}
  0\\ 1
\end{bmatrix},\\
G_{11}=S_\eta, \quad G_{21}=\begin{bmatrix}
  1 \\ 1
\end{bmatrix}, \quad  K_{x1}=-1, \\
K_{\zeta 1}=[-500,\,\,-500], 
G_{12}=\begin{bmatrix}
  0 & 1 \\
  -w^2 & 0
\end{bmatrix}, \quad G_{22}=\begin{bmatrix}  0 \\ 1
\end{bmatrix}, \\ K_{x2}=-2, \quad K_{\zeta 2}=[0, \,\, -500]. 
\end{gather*}
It can be verified that assumption A5) is satisfied for both source nodes under the controller gains defined above. All the conditions in Theorem \ref{th05} are satisfied, therefore,
there is $\epsilon$ for the above controller to make the electrical network realize the output cooperation.

\begin{figure}
  \centering
 \includegraphics[width=0.45\textwidth]{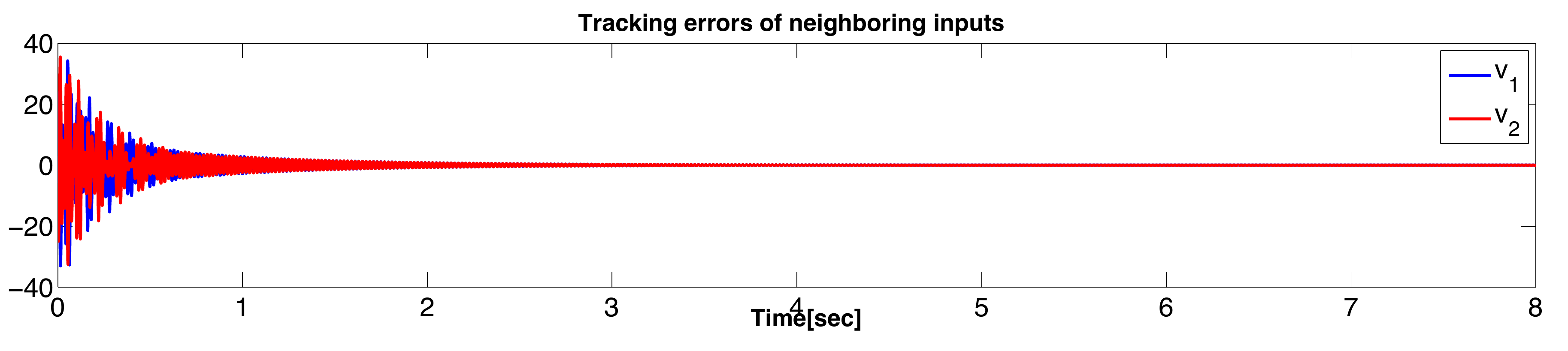}\\
  \includegraphics[width=0.45\textwidth]{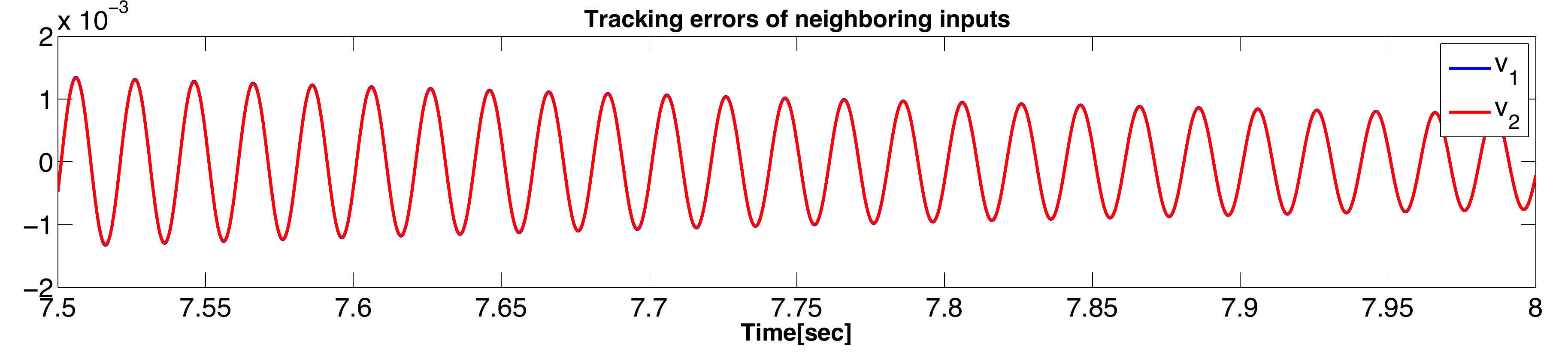}\\
    \includegraphics[width=0.45\textwidth]{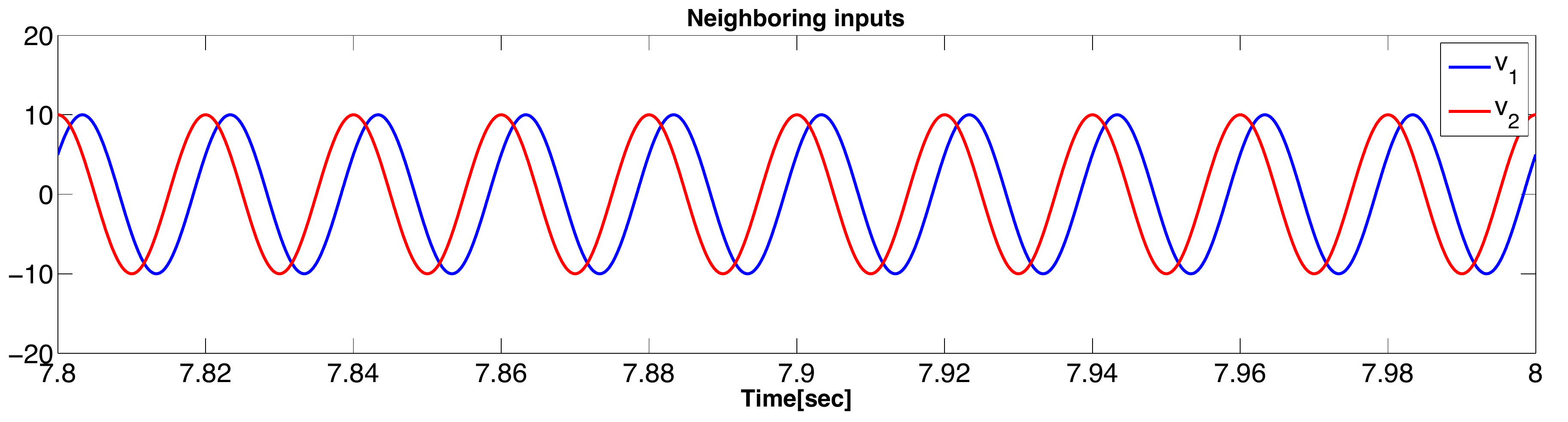}\\
  \caption{The trajectories of the neighboring inputs of both nodes and their tracking errors; the top: the trajectories of tracking errors; the middle: an enlargement of the top; the bottom: an enlargement of the trajectories of the neighboring inputs.}\label{fig05}
\end{figure}

Simulation results with $\epsilon=20$ made in the Matlab environment are shown in Fig. \ref{fig05}. The electrical network is built by making use of the SimPowerSystems Toolbox. After the transition time, each generator adjusts its output to make the neighboring input of itself converge to the desired one. As shown in the middle figure of Fig. \ref{fig05}, the tracking error is in the order of $10^{-3}$ and is still decreasing. A less $\epsilon$ leads to a longer transition phase, but the error will always ultimately decay to zero.

\section{Conclusion} \label{sec07}
A new class of multi-agent network systems was presented, where the nodes are not directly coupled but indirectly coupled by dynamic systems, called dynamic edges. The node dynamics can be directly controlled and are influenced by the neighboring input which is a weighted sum of the edge outputs; while the edge dynamics can not be directly controlled due to its input being the node outputs. Distributed controllers designed by a combination of feedback passivity theories and the internal model principle were presented for output synchronization, output cooperation and master-slave output cooperation, respectively. A simulation example of cooperative current control of an electrical network illustrates the efficacy of the analytic results.

The developments were based on the exact matrices, but they can be extended to the uncertain case in that the tools of both passification and internal model principle are good at coping with the uncertain systems. Although the network is heterogeneous, the controller requires some common parameters, e.g., $B_\eta$, $C_\eta$ and $\epsilon$. How to relax such a common requirement and how to get the upper bound of $\epsilon$ are the goals of our future researches.

\appendix
\section{A preliminary lemma}
\begin{lemma} \label{le01}
Given a matrix $W = \begin{bmatrix}W_1 & W_2+ W_5 \\
W_3 & W_4 \end{bmatrix}$. If there are symmetric positive definite matrices $P_w$ and $Q_w$ such that
\begin{subequations} 
\begin{align}
   P_wW_1+W_1^TP_w\le 0,\\
   Q_wW_4+W_4^TQ_w<0, \label{A02b} \\ 
   P_w W_2=-W_3^TQ_w,
\end{align}
\end{subequations}
and $W_1$ is Hurwitz, then there is a constant $\bar{\epsilon}$ such that for all $\|W_5\|<\bar{\epsilon}$ there is a block diagonal symmetric positive definite matrix $\bar{P} = \mathrm{diag}(\bar{P}_w, Q_w)$ for some $\bar{P}_w>0$ satisfying
\begin{equation}  \label{A03}
  \bar{P}W+W^T\bar{P}<0.
\end{equation}
\end{lemma}

\begin{IEEEproof}
The proof is constructive. Since $W_1$ is Hurwitz, there is a symmetric positive definite matrix $P_r$ satisfying the following Lyapunov equation,
\begin{equation}
  P_rW_1+W_1^TP_r = -I.
\end{equation}
By \eqref{A02b}, there is positive scalar $\epsilon_1$ such that $Q_wW_4+W_4^TQ_w<-\epsilon_1$. Take $\bar{P}_w=P_w+\bar{\epsilon}P_r$ to obtain
\begin{multline}
 \bar{P}W+W^T\bar{P} = \\ \begin{bmatrix}
   -\bar{\epsilon}I+P_wW_1+W_1^TP_w & \star \\
   (P_w W_5 + \bar\epsilon P_r(W_2+W_5))^T & Q_wW_4+W_4^TQ_w
 \end{bmatrix}\\
 <\begin{bmatrix}
   -\bar{\epsilon}I & \star \\
   (P_w W_5 + \bar\epsilon P_r(W_2+W_5))^T & -\epsilon_1I
 \end{bmatrix}
\end{multline}
where $\star$ denotes the symmetric part.
By Finsler lemma \cite{BoydBOOK1994}, it can be verified that the above matrix is symmetric negative definite, if 
\begin{equation*}
-\epsilon_1 I+\bar{\epsilon}^{-1}\|P_w W_5 + \bar\epsilon P_r(W_2+W_5)\|^2<0, 
\end{equation*}
which can be further enlarged by
\begin{equation} \label{RA05}
a_r^2\bar{\epsilon}^3+2a_wa_r\bar{\epsilon}^2+a_w^2\bar{\epsilon}-\epsilon_1<0,
\end{equation}
where $a_r=\|P_r\|$ and $a_w=\|P_w\|+\|P_rW_2\|$. Selecting $\bar{\epsilon}<\min\{1,\frac{\epsilon_1}{(a_r+a_w)^2}\}$ makes \eqref{RA05} true to complete the proof. 

Actually, if $a_r$ and $a_w$ are such that $(1-a_r)\sqrt[3]{\epsilon_1}>a_w$, then all $\bar{\epsilon}<\sqrt[3]{\epsilon_1}$ satisfies \eqref{RA05}. This means that $\bar{\epsilon}$ is not always restricted to be a small value since if $\epsilon_1$ is a large value, then $\bar{\epsilon}$ can be as well a high value. 
\end{IEEEproof}

\section{Proof of Theorems}
\subsection{Proof of Theorem \ref{th00}}
\begin{IEEEproof}
As shown in \cite{AndrievskyCDC1996}, the hyper minimum phase condition implies that there is a SPD matrix $P_s$ and a matrix $K_{xi}$ such that 
$(A_i+B_iK_{xi})P_s+P_s(A_i+B_iK_{xi})^T<0$ and $C_iP_s=D_i^T=B_i^T$.
Meanwhile, the incorporation of p-copy internal model of $S_\eta$ and assumption A2) imply that there is a SPD matrix $P_g$ such that $G_{1i}P_g+P_gG_{1i}^T\le 0$. Design $K_{\zeta i}$ by $K_{\zeta i}=-G_{2i}^TP_g^{-1}$. Let $\hat{P}=\mathrm{diag}(P_s,P_g)$. It can be verified that 
\begin{equation} \label{RA01}
\hat{A}\hat{P}+\hat{P}\hat{A}^T\le 0, \quad \hat{C}_i\hat{P}=\hat{D}_i,
\end{equation}
which means that system $(\hat{A}_i, \hat{D}_i,\hat{C}_i)$ is passive. 
Let $G_{2i}$ be such that not only $(G_{1i}, G_{2i})$ is controllable but also $(G_{i1}, G_{2i}^TP_g^{-1})$ is observable, which together with \eqref{RA01} implies that $\hat{A}$
is Hurwitz. 
\end{IEEEproof}

\subsection{Proof of Theorem \ref{th01}}
\begin{IEEEproof}
The closed-loop system has the following compact form
\begin{equation} \label{A05}
   \dot{\hat{X}} = \hat{\mathcal{A}}\hat{X} + \hat{\mathcal{D}}_\eta\eta, \quad \dot{\eta}=(I_N\times S_\eta)\eta,
\end{equation}
where $\hat{X}=[\hat{x}^T_1,\cdots,\hat{x}_N^T,z^T_1,\cdots,z_M^T]^T$, $\eta = [\eta_1^T, \eta_2^T,\cdots, \newline \eta_N^T]^T$,
\begin{gather*}
\hat{\mathcal{A}} = \begin{bmatrix}
  \overline{AN} & -\overline{HDG} \\
  \overline{H^TFC} & \overline{EM}
\end{bmatrix},
\overline{AN} = \begin{bmatrix}
  \hat{A}_1 & \cdots & 0 \\
   \vdots & \ddots & \vdots &  \\
  0 & \cdots & \hat{A}_N
  \end{bmatrix},  \\ \overline{HDG} = \begin{bmatrix}
    h_{11}\hat{D}_1G_1 & \cdots & h_{1M}\hat{D}_1G_M\\
  \vdots & \vdots & \vdots \\
  h_{N1}\hat{D}_NG_1 & \cdots & h_{NM}\hat{D}_NG_M
  \end{bmatrix},\nonumber\\
 \overline{H^TFC} = \begin{bmatrix}
    h_{11}F_1\hat{C}_1 & \cdots & h_{N1}F_1\hat{C}_N  \\
  \vdots & \vdots & \vdots \\
  h_{1M}F_M\hat{C}_1 & \cdots & h_{NM}F_M\hat{C}_N
  \end{bmatrix}, \\ 
  \overline{EM} = \begin{bmatrix}
  {E}_1 & \cdots & 0 \\
   \vdots & \ddots & \vdots &  \\
  0 & \cdots & {E}_M
  \end{bmatrix}, \hat{\mathcal{D}}_\eta =  \begin{bmatrix}
  \hat{D}_{\eta 1} & \cdots & 0 \\
   \vdots & \ddots & \vdots &  \\
  0 & \cdots & \hat{D}_{\eta N}  \end{bmatrix}. 
 \end{gather*}
According to Lemma \ref{le01a}, in order for $e_i\rightarrow 0$, it suffices to show $\hat{\mathcal{A}}$ is Hurwitz.

With assumption A3), there exists a positive definite matrix $Q_j$ such that
\begin{equation}  \label{A07}
  Q_jE_j+E_j^TQ_j<0, \quad Q_jF_j=G_j^T, \quad j=1,2,\cdots,M,
\end{equation}
and with assumption A5), there exists a positive definite matrix $\hat{P}_i$ such that
\begin{equation}  \label{A08}
  \hat{P}_i\hat{A}_i+\hat{A}_i^T\hat{P}_i\le 0, \quad \hat{P}_i\hat{D}_i=\hat{C}_i^T,\quad i=1,2,\cdots,N,
\end{equation}
Define $P=\mathrm{diag}(\hat{P}_1,\cdots,\hat{P}_N)$ and $Q = \mathrm{diag}(Q_1,\cdots,\break Q_M)$. It can be seen that $P\overline{HDG}=(Q\overline{H^TFC})^T$. Therefore by Lemma \ref{le01}, $\hat{\mathcal{A}}$ is Hurwitz, and subsequently all $z_i$, $x_i$ are bounded for bounded external inputs $\eta_i$. Moreover, $e_i\rightarrow 0$ for all $i=1,2,\cdots,N$ by Lemma \ref{le01a} due to the incorporation of a $p$-copy internal model in controller \eqref{III03}.
\end{IEEEproof}

\subsection{Proof of Theorem \ref{th02}}
\begin{IEEEproof}
The closed-loop system with controller \eqref{III07} has the form of
\begin{equation} \label{A10}
  \left\{
  \begin{split}
    &\dot{\hat{x}}_i = \hat{A}_i \hat{x}_i - \hat{D}_i\sum_{j=1}^M h_{ij} w_j + \hat{D}_{\eta i}\eta_i,\quad i=1,\cdots,N, \\
    &\dot{z}_j = E_jz_j + F_j \sum_{i=1}^N h_{ij}y_i, \,\, j=1,\cdots, M, \\
    &\dot{\eta}=(I_N\otimes S_\eta)\eta - \epsilon (H\otimes B_\eta) w,
  \end{split}
  \right.
\end{equation}
Consider error vectors $e_{xi}=\hat{x}_i-\Pi_i \eta_i$. From \eqref{III05}, it follows that
\begin{equation} \label{A11}
\dot{e}_{xi}=\hat{A}_ie_{xi}-\hat{D}_i\sum_{j=1}^M h_{ij} w_j + \epsilon \Pi_iB_\eta \sum_{j=1}^M h_{ij} w_j, 
\end{equation}
for all $i=1,\cdots,N$. 
Noticing \eqref{III06}, the dynamics of $z_j$ can be written by
\begin{multline} \label{A12}
 \dot{z}_j = E_jz_j + F_j \sum_{i=1}^N h_{ij}\hat{C}_ie_{xi}+F_j \sum_{i=1}^N h_{ij} Q_\eta \eta_i, 
 \end{multline}
for all $j=1,\cdots,M$. Recall matrix $T\in\real^{(N-1)\times N}$ satisfying $T\underline{\mathbf{1}}=0$ and $TT^T=I_{N-1}$, which has been given after assumption A4), and introduce a coordinate transformation by $\tilde{\eta}_i=(T_i\otimes I_q)\eta$, $i=1,\cdots,N-1$, and $\tilde{\eta}_N=\sum_{i=1}^N\eta_i$, where $T_i$ denotes the $i$th row of $T$. The third formula in \eqref{A10} can be transformed into
\begin{equation} \label{A13}
  \left\{
  \begin{split}
    &\dot{\tilde{\eta}}_i=S_\eta \tilde{\eta}_i -\epsilon B_\eta \sum_{j=1}^M \bar{h}_{ij} w_j, \,\, i=1,\cdots\!,N-1. \\
    &\dot{\tilde{\eta}}_N = S_\eta\tilde{\eta}_N
  \end{split}
  \right.
\end{equation}
where $\bar{h}_{ij}$ is the $i$th row and $j$th column element of matrix $\bar{H}$. Accordingly, equation \eqref{A12} can be rewritten by
\begin{multline} \label{A14}
 \dot{z}_j = E_jz_j + F_j \sum_{i=1}^N h_{ij}\hat{C}_ie_{xi}+F_j \sum_{i=1}^{N-1} \bar{h}_{ij} Q_\eta \tilde{\eta}_i.
\end{multline}
Since $\tilde{\eta}_N$ freely runs, it suffices to consider the dynamics of $e_{xi}$, $z_j$ and $\tilde{\eta}_k$, with $i=1,\cdots,N$, $j=1,\cdots,M$ and $k=1,\cdots,N-1$, respectively.
Define a stacked vector by
$\bar{X}=[e^T_{x1}, \cdots, e_{xN}^T, $ $z_1^T, \cdots, z_M^T, \tilde{\eta}_1^T, \cdots, \tilde{\eta}_{N-1}^T]^T$, whose dynamics is governed by
\begin{equation} \label{A15}
 \dot{\bar{X}}=\bar{\mathcal{A}}\bar{X}
\end{equation}
where
\begin{gather*}
\bar{\mathcal{A}} = \begin{bmatrix}
\overline{AN} & -\overline{HDG}+\epsilon \overline{\Pi B} & 0 \\
\overline{H^TFC} & \overline{EM} & \overline{\bar{H}^TF}(I_{N-1}\otimes Q_\eta) \\
0 & -\epsilon(I_{N-1}\otimes B_\eta)\overline{\bar{H}G} & I_{N-1}\otimes S_\eta
\end{bmatrix}
\\
\overline{\bar{H}^TF}=\begin{bmatrix}
 \bar{h}_{11}F_1& \cdots & \bar{h}_{(N-1)1}F_1\\
   \vdots & \vdots & \vdots\\
\bar{h}_{1M}F_M & \cdots & \bar{h}_{(N-1)M}F_M
\end{bmatrix},  \nonumber \\
\overline{\bar{H}G}=\begin{bmatrix}
\bar{h}_{11} G_1 & \cdots & \bar{h}_{1M} G_M \\
 \vdots & \vdots & \vdots \\
 \bar{h}_{(N-1)1} G_1 & \cdots & \bar{h}_{(N-1)M} G_M
\end{bmatrix}, \\ \overline{\Pi B} =  \begin{bmatrix}
{h}_{11}\Pi_1 B_\eta G_1 & \cdots & {h}_{1M}\Pi_1B_\eta G_M \\
 \vdots & \vdots & \vdots \\
 {h}_{N1}\Pi_N B_\eta G_1 & \cdots & \bar{h}_{NM}\Pi_N B_\eta G_M
\end{bmatrix}
\end{gather*}
With \eqref{A07} and \eqref{A08} and using Lemma \ref{le01}, there is upper bound $\bar{\epsilon}$ such that for all $0<\epsilon<\bar{\epsilon}$,  such that
\begin{multline}
  \mathcal{P}\begin{bmatrix}
\overline{AN} & -\overline{HDG}+\epsilon \overline{\Pi B} \\
\overline{H^TFC} & \overline{EM}
\end{bmatrix}+ \\ \begin{bmatrix}
\overline{AN} & -\overline{HDG}+\epsilon \overline{\Pi B} \\
\overline{H^TFC} & \overline{EM}
\end{bmatrix}^T\mathcal{P}<-\epsilon_1 I
\end{multline}
for some positive scalar $\epsilon_1$, where $\mathcal{P}=\mathrm{diag}(\tilde{P}_1,\cdots,\tilde{P}_N,\break Q_1,\cdots,Q_M)$ is a block diagonal symmetric positive matrix with $Q_iF_i=G_i^T$, for all $i=1,\cdots, M$.

Consider the Lyapunov function $\hat{V} = \sum_{i=1}^N e_{xi}^T \tilde{P}_ie_{x i}+ \sum_{j=1}^M z_j^T Q_jz_j
+ \epsilon^{-1}\sum_{i=1}^{N-1} \tilde{\eta}_i^TP_\eta \tilde{\eta}_i$. Its time derivative along system \eqref{A15} yields
\begin{multline} \label{A18}
  \dot{\hat{V}} \le -\epsilon_1\sum_{i=1}^N {e}_{xi}^T{e}_{xi}
   -\epsilon_1\sum_{j=1}^M {z}_i^T{z}_i \\ + \epsilon^{-1} \sum_{i=1}^{N-1}\tilde{\eta}^T_i(P_\eta S_\eta+S_\eta^T P_\eta)\tilde{\eta}_i\le 0,
\end{multline}
by which, the invariant set of $\mathcal{D}=\{\bar{X}:\dot{\hat{V}}=0\}$ is a subset of 
$\{\bar{X}:{e}_{xi}=0,z_j=0, \forall i,j\}$. According to the Lasalle's theorem,  as $t\rightarrow \infty$, $\bar{X}(t)\in \mathcal{D} $, so $e_{xi}\rightarrow 0$ and $z_j\rightarrow 0$ for all $i,j$. The latter means that
$s_j\rightarrow 0$ due to $F_j$ of full column rank, which in a compact form is
\begin{equation} \label{A19} 
(H^T\otimes I_p)Y\rightarrow 0
\end{equation}
where $Y=[y_1^T,y_2^T,\cdots,y_N^T]^T$. With assumption A4), the null space of $H^T$ has dimension $1$ and is spanned by vector $\underline{\mathbf{1}}$, therefore,
\begin{equation}
  \label{A20} 
y_i\rightarrow y_k, \quad \forall \,  i,\,k,
\end{equation}
which, together with $e_{xi}\rightarrow 0$, means that
\begin{equation}\label{A21} 
Q_\eta\eta_i-Q_\eta\eta_k\rightarrow 0,    \quad \forall \,  i,\,k.
\end{equation}
Notice that the sum of all $\eta_i$s, namely $\tilde{\eta}_N$, has the following dynamics
\begin{equation} \label{A22} 
\dot{\tilde{\eta}}_N=S_\eta \tilde{\eta}_N, \quad \mathrm{with}\quad \tilde{\eta}_N(0)=\sum_{i=1}^N \eta_i(0).
\end{equation}
Combining \eqref{A21} and \eqref{A22} yields
\begin{equation}
  y_i\rightarrow Q_\eta \eta_i\rightarrow \frac{1}{N}Q_\eta \tilde{\eta}_N,
\end{equation}
which is just \eqref{III08}.
\end{IEEEproof}

\subsection{Proof of Theorem \ref{th03}}
\begin{IEEEproof}
According to Lemma \ref{le01a}, it is enough to show the internal stability of system \eqref{IV09a} in which state variable $Z$ is driven by exosystem $\tilde{\nu}$. Consider the Lyapunov function $V_1 = \sum_{j=1}^Mz_j^TQ_jz_j+\epsilon^{-1}\sum_{i=1}^{N-1}\tilde{\eta}_i^T(I_p\otimes P_\eta)\tilde{\eta}_i$, whose time derivative along \eqref{IV09a} with $\nu=0$ is given by
\begin{equation} \label{A24}
  \dot{V}_1 \le \sum_{j=1}^Mz_j^T(Q_jE_j+E_j^TQ_j)z_j\le 0.
\end{equation}
By Lasalle's lemma, $z_j\rightarrow 0$, which together with $F_j$ being of full column rank, means that $\bar{H}^T(I_{N-1}\otimes G_Q)\tilde{\eta}\rightarrow 0$. Furthermore, because of $\bar{H}^T$ being of full column rank, $G_Q\tilde{\eta}_i\rightarrow 0$. Noting that matrix pair $(G_S, G_Q)$ is observable, one has $\tilde{\eta}_i\rightarrow 0$. This means the system \eqref{IV09a} is internal stable. Notice matrix pair $(G_S, G_B)$ is a $p$-copy model of $S_\eta$, hence according to Lemma 1
\begin{equation} \label{A25}
-\sum_{j=1}^M\bar{h}_{ij} G_jz_j-Q_v\tilde{\nu}_i\rightarrow 0,
    \end{equation}
which means that \eqref{IV10} holds.

Notice that \eqref{A25} is equivalent to $(T\otimes I)(v-\bar{Q}_v\nu)\rightarrow 0$ with $\bar{Q}_v=(I_N\otimes Q_v)$. Therefore
\begin{equation}
  v-\bar{Q}\nu\rightarrow \underline{\mathbf{1}}\otimes \nu_0
\end{equation}
for some $\nu_0\in\real^p$. By the fact $(\underline{\mathbf{1}}^T\otimes I_p)v=0$, it follows that
\begin{equation} \label{A27}
  -Q_v\sum_{i=1}^N\nu_i = -Q_v\tilde{\nu}_N\rightarrow N\nu_0,
\end{equation}
that is, $\nu_0\rightarrow -\frac{1}{N}Q_v\tilde{\nu}_N$,
by which \eqref{IV11} is obtained.
\end{IEEEproof}

\subsection{Proof of Theorem \ref{th04}}
\begin{IEEEproof}
With controller \eqref{IV02}, the closed-loop system has the form of
\begin{equation} \label{A28}
  \left\{
  \begin{split}
    &\dot{\hat{x}}_i = \hat{A}_i \hat{x}_i - \hat{D}_i\sum_{j=1}^M h_{ij} G_jz_j + \hat{D}_{\bar{\eta} i}\bar{\eta}_i, \,\, i=1,2,\cdots,N, \\
    &\dot{z}_j = E_jz_j + F_j \sum_{i=1}^N h_{ij}y_i, \quad j=1,2,\cdots, M, \\
    &\dot{\bar{\eta}}_i=G_S \bar{\eta}_i - \epsilon G_B \sum_{j=1}^M h_{ij} G_jz_j -\epsilon G_B Q_v\nu_i, \\[-1em] &\quad \quad \quad \quad \quad \quad \quad \quad \quad \quad \quad \quad \hspace{3em}i=1,2,\cdots,N.\\
    &\dot{\nu}_i = S_\eta\nu_i, \quad i=1,2,\cdots,N.
  \end{split}
  \right.
\end{equation}
With the same notations as those in \eqref{IV08} and \eqref{IV09}, the last three equations in \eqref{A28} can be transformed into,
\begin{subequations}
  \label{A29}
  \begin{align}
    & \dot{z}=\overline{EM}z + \overline{\bar{H}^TF}(I_{N-1}\otimes G_Q)\tilde{\eta} \nonumber\\
    &\quad \quad\quad \quad \quad +\begin{bmatrix}
      F_1\sum_{i=1}^N h_{i1}(y_i-G_Q\bar{\eta}_i) \\ \vdots \\ F_M\sum_{i=1}^N h_{iM}(y_i-G_{Q}\bar{\eta}_i)
    \end{bmatrix}, \\
    & \dot{\tilde{\eta}} = (I_{N-1}\otimes G_S)\tilde{\eta} \nonumber \\
   &\quad \quad  -\epsilon (I_{N-1}\otimes G_B)(\overline{\bar{H}G}z+(I_{N-1}\otimes Q_v)\tilde{\nu}),\\
    & \dot{\tilde{\eta}}_N = G_S\tilde{\eta}_N, \quad \tilde{\nu}_N\equiv 0,
  \end{align}
\end{subequations}
where the last equation is from assumption A6). Consider error vectors $\bar{e}_{xi}=\hat{x}_i-\bar{\Pi}_{1i}\bar{\eta}_i-\bar{\Pi}_{2i}\nu_i$, where $\bar{\Pi}_{1i}$ and $\bar{\Pi}_{2i}$ are the solution of \eqref{IV06} and \eqref{IV07}. Its dynamics has the form of,
\begin{multline}
  \dot{\bar{e}}_{xi}=\hat{A}_i\bar{e}_{xi}-\hat{D}_i\left(\sum_{j=1}^M h_{ij}G_jz_j + Q_v \nu_i \right)\\ +
  \epsilon \bar{\Pi}_{1i}G_B\left(\sum_{j=1}^M h_{ij}G_jz_j+Q_v\nu_i \right).
\end{multline}

Define
\begin{gather*}
  \overline{HG}=\begin{bmatrix}
    h_{11}G_1 & \cdots  & h_{1M}G_M \\
    \vdots &    \vdots  & \vdots \\
   h_{N1}G_1 &  \cdots  & h_{NM}G_M
  \end{bmatrix},\\
  \overline{\Pi_1G_B} = \begin{bmatrix}
  h_{11}\Pi_{11}G_BG_1  &   \cdots  & h_{1M}\Pi_{11}G_BG_M \\
  \vdots                &   \vdots  & \vdots\\
  h_{N1}\Pi_{1N}G_BG_1  &   \cdots  & h_{NM}\Pi_{1N}G_BG_M
  \end{bmatrix}.
\end{gather*}

Consider errors $\bar{e}_x = [\bar{e}_{x1}^T, \cdots, \bar{e}_{xN}^T]^T$, $\bar{e}_z = z-\tilde{\Pi}_z\tilde{\nu}$ and $\bar{e}_\eta = \tilde{\eta}-\tilde{\Pi}_\eta\tilde{\nu}$, where $\tilde{\Pi}_z$ and $\tilde{\Pi}_\eta$ are the solution of \eqref{IV10}. Noting that $y_i-G_Q\bar{\eta}_i=\hat{C}_i\bar{e}_{xi}$ due to \eqref{IV07} and
\begin{multline} \label{A32}
  \begin{bmatrix}
    T\otimes I_p \\
    \underline{\mathbf{1}}^T\otimes I_p
  \end{bmatrix}\big(\overline{HG}\tilde{\Pi}_z\tilde{\nu} +(I_N\otimes Q_v)\nu\big) \\ = \begin{bmatrix}
    \overline{\bar{H}G}\tilde{\Pi}_z \tilde{\nu} + (I_{N-1}\otimes Q_v)\tilde{\nu} \\
     Q_v\tilde{\nu}_N
  \end{bmatrix} = 0,
\end{multline}
their dynamics is governed by
\begin{equation}
\left\{
\begin{split}
  &\dot{\bar{e}}_x = \overline{AN}\bar{e}_x - \overline{HDG}\bar{e}_z + \epsilon \overline{\Pi_1G_B}\bar{e}_z \\
  &\dot{\bar{e}}_z = \overline{EM}\bar{e}_z + \overline{H^TFC}\bar{e}_x + \overline{\bar{H}^TF}(I_{N-1}\otimes G_Q)\bar{e}_\eta \\
  & \dot{\bar{e}}_\eta = -\epsilon(I_{N-1}\otimes G_B)\overline{\bar{H}G}\bar{e}_z + (I_{N-1}\otimes G_s)\bar{e}_\eta
\end{split}\right..
\end{equation}
Making use of Lemma \ref{le01}, it can be verified that the above system is Hurwitz along with the line of the proof of Theorem \ref{th02}. From \eqref{A32}, it follows that
\begin{equation}
  \overline{HG}\tilde{\Pi}_z\tilde{\nu}=-(I_N\otimes Q_v)\nu,
\end{equation}
which, together with $\bar{e}_z\rightarrow 0$, implies that $v_i\rightarrow Q_v\nu_i$ for all $i$, and thus the output cooperation is realized.

On the other hand, from $\tilde{\eta}_N=\sum_{i=1}^N\bar{\eta}_i$ and $\bar{e}_{xi}\rightarrow 0$, it follows that \eqref{IV12} holds.
\end{IEEEproof}

\subsection{Proof of Theorem \ref{th05}}
\begin{IEEEproof}
Firstly, we consider the edge dynamics by replacing $y_i$ by $G_Q\bar{\eta}_i$ for $i=1,\cdots,l$ and by $Q_\eta\eta_i$ for $i=l+1,\cdots,N$.
\begin{equation}
\left\{
  \begin{split}
    &\dot{z}_j = E_jz_j + F_j\sum_{i=1}^{l}h_{ij}G_Q\bar{\eta}_i\\
    &\quad\quad \quad \quad \quad  + F_j\sum_{k=l+1}^Nh_{kj}Q_\eta\eta_k, \quad j=1,\cdots,M, \\
    & \bar{\eta}_i=G_S\bar{\eta}_i + \epsilon G_B(-\sum_{j=1}^M h_{ij}G_jz_j-Q_v \nu_i),
    \\ &\quad\quad \quad \quad \quad\quad \quad \quad \quad\quad \quad \quad   i=1,\cdots, l,\\
    & \dot{\nu}_i=S_\eta \nu_i, \quad i=1,\cdots,l,\\
    & \dot{\eta}_k=S_\eta \eta_k, \quad k=l+1,\cdots,N.
  \end{split}\right.
\end{equation}
Notice that $H^T\begin{bmatrix}
  I_l \\ 0
\end{bmatrix}$ is of full column rank, hence by the similar way in the above theorems it can be derived that the system consisting of $z_j$, $j=1,\cdots,M$ and $\bar{\eta}_i$, $i=1,\cdots,l$, is internal stable. Therefore, there are unique matrices $\Pi^\nu_{z j}$, $\Pi_{z j}^\eta$, $\Pi_{\bar{\eta} i}^\nu$ and $\Pi_{\bar{\eta} i}^\eta$ satisfying the following matrix equation
\begin{multline}
 \Pi_{zj}^\nu (I_l\otimes S_\eta)= E_j\Pi_{zj}^\nu + F_j\sum_{i=1}^l h_{ij}G_Q\Pi_{\bar{\eta}i}^\nu, \\ j = 1, \cdots,M, 
\end{multline}
\begin{multline}
\Pi_{z j}^\eta(I_{N-l}\otimes S_\eta) = E_j\Pi_{z j}^\eta + F_j\sum_{i=1}^l h_{ij}G_Q \Pi_{\bar{\eta} i}^\eta \\ +
F_j\sum_{k=l+1}^N h_{kj}Q_\eta(\xi^{k-l}_{N-l}\otimes I_q), \,\, j=1,\cdots,M,
\end{multline}
\begin{multline}
\Pi_{\bar{\eta} i}^\eta (I_{N-l}\otimes S_\eta) = G_S\Pi_{\bar{\eta} i}^\eta, \quad \Pi_{\bar{\eta} i}^\nu(I_{l}\otimes S_\eta)=G_S\Pi_{\bar{\eta} i}^{\nu},\\  i=1,\cdots,l,
\end{multline}
\begin{multline}
-\sum_{j=1}^M h_{ij}G_j\Pi_{z j}^{\nu}=\xi_l^i\otimes Q_v,\quad \sum_{j=1}^M h_{ij}G_j\Pi_{z j}^{\eta}=0, \\ \quad i=1,\cdots, l, \label{A39}
\end{multline}
such that
\begin{equation} \label{A40}
\begin{aligned}
& z_j\rightarrow \Pi_{z j}^\nu \hat{\nu}+\Pi_{z j}^\eta\hat{\eta}, &j& =1,\cdots,M,\\
&\bar{\eta}_i\rightarrow \Pi_{\bar{\eta} i}^\nu \hat{\nu}+\Pi_{\bar{\eta} i}^{\eta}\hat{\eta},
&i&=1,\cdots,l,
\end{aligned}
\end{equation}
where $\hat{\nu}=[\nu^T_1,\cdots,\nu^T_l]^T$, $\hat{\eta}=[\eta_{l+1}^T,\cdots, \eta_{N}^T]^T$, $I_l$ is a unit matrix with order $l$, and henceforth notation $\xi_n^m$ denotes the $n$-order row vectors with all elements being 0 but the $m$th element being $1$,  for example, $\xi_l^i=[0,\cdots,\underset{i}{1},\cdots, 0]\in\real^l$.

Secondly, consider the node dynamics. The closed-loop node dynamic system has the form of
\begin{equation}
  \left\{
  \begin{array}
    {l}
    \dot{\hat{x}}_i = \hat{A}_i\hat{x}_i + \hat{D}_i v_i + \hat{D}_{\bar{\eta}_i}\bar{\eta}_i, \quad i=1,\cdots, l,\\
    \dot{\hat{x}}_i = \hat{A}_i\hat{x}_i + \hat{D}_i v_i + \hat{D}_{\eta_i}\eta_i, \quad i=l+1,\cdots, N.\\
  \end{array}
  \right.
\end{equation}
Replacing $v_i$ by $-\sum_{j=1}^Mh_{ij}G_j(\Pi_{z j}^\nu \hat{\nu}+\Pi_{z j}^\eta \hat{\eta})$, for all $i=1,\cdots, N$, yields
\begin{equation}
\left\{\begin{array}{l}
  \dot{\hat{x}}_i=\hat{A}_i\hat{x}_i+ \hat{D}_iM_i \hat{\nu} + \hat{D}_{\bar{\eta}i} \bar{\eta}_i, \quad i=1,\cdots,l\\
  \dot{\hat{x}}_i=\hat{A}_i\hat{x}_i+ \hat{D}_iM_i \hat{\nu} + \hat{D}_i N_i \hat{\eta}+\hat{D}_{\eta i}(\xi_{N-l}^{i-l}\otimes I_q)\hat{\eta}, \\
  \quad \quad \quad \quad \quad \quad \quad \quad \quad  i=l+1,\cdots,N,
  \end{array}
  \right.
\end{equation}
where $M_i=-\sum\limits_{j=1}^M h_{ij}G_j\Pi_{z j}^{\nu}$ and $N_i=-\sum\limits_{j=1}^M h_{ij}G_j\Pi_{z j}^{\eta}$, $i=1,\cdots,N$ (Please note $N_i=0$ for $i=1,\cdots,l$).  Since $\hat{A}_i$ is Hurwitz, there are unique matrices $\Pi^\nu_{fi}$ and $\Pi_{fi}^{\bar{\eta}}$, $i=1,\cdots, l$, satisfying
\begin{equation} \label{A42}
\left\{\begin{split}
  &\Pi_{fi}^{\nu}(I_{l}\otimes S_\eta)=\hat{A}_i\Pi_{fi}^\nu+\hat{D}_iM_i, \\ &\Pi_{fi}^{\bar{\eta}}G_s=\hat{A}_i\Pi_{fi}^{\bar{\eta}}+\hat{D}_{\bar{\eta}i}, \\
  & \hat{C}_i\Pi_{fi}^\nu = 0, \quad \hat{C}_i\Pi_{fi}^{\bar{\eta}} = G_Q,
\end{split}\right.
\end{equation}
for the first $l$ slave nodes, and there are unique matrices $\Pi^\nu_{li}$ and $\Pi_{li}^\eta$, $i=l+1,\cdots,N$, satisfying
\begin{equation} \label{A44}
\left\{
\begin{split}
  &\Pi_{li}^{\nu}(I_{l}\otimes S_\eta)=\hat{A}_i\Pi_{li}^\nu+\hat{D}_iM_i,\\
  &
  \Pi_{li}^{\eta}(I_{N-l}\otimes S_\eta)=\hat{A}_i\Pi_{li}^\eta+\hat{D}_iN_i+ \xi_{N-l}^{i-l}\otimes \hat{D}_{\eta i} \\
  & \hat{C}_i\Pi^\nu_{li}=0, \quad \hat{C}_i\Pi_{li}^\eta=\xi_{N-l}^{i-l}\otimes Q_\eta \\
  \end{split} \right.
\end{equation}
for the remainder $N-l$ master nodes. The last equation of \eqref{A42} and \eqref{A44} comes from the incorporation of a $p$-copy internal model in \eqref{V01a} and \eqref{V01b}, respectively.

Finally, we are ready to consider the whole closed-loop system under controller \eqref{V01}, which has the form of
\begin{equation*}
\left\{\begin{array}{ll}
  \dot{\hat{x}}_i = \hat{A}_i \hat{x}_i + \hat{D}_i v_i + \hat{D}_{\bar{\eta}_i}\bar{\eta}_i, & i=1,\cdots,l, \\
  \dot{\bar{\eta}}_i = G_S\bar{\eta}_i + \epsilon G_B(v_i-Q_v\nu_i), &i=1,\cdots,l, \\
  \dot{\hat{x}}_i = \hat{A}_i \hat{x}_i + \hat{D}_i v_i + \hat{D}_{{\eta}_i}{\eta}_i, & i=l+1,\cdots,N, \\
  \dot{z}_j = E_jz_j + F_j\sum_{i=1}^N h_{ij}y_i, & j=1,\cdots,M, \\
  \dot{\hat{\nu}}=(I_l\otimes S_\eta)\hat{\nu}, & \hat{\nu}=[\nu^T_1,\cdots,\nu^T_l]^T,\\
  \dot{\hat{\eta}} = (I_{N-l}\otimes S_\eta)\hat{\eta},     &  \hat{\eta}=[\eta_{l+1}^T,\cdots, \eta_{N}^T]^T.
  \end{array}
\right.
\end{equation*}
Consider the following error vectors, 
\begin{align}
&\hat{e}_{xi}=\left\{\begin{array}
  {ll} \hat{x}_i-\Pi_{fi}^\nu\hat{\nu}-\Pi_{fi}^{\bar{\eta}}\bar{\eta}_i,& i=1,\cdots, l\\
  \hat{x}_i-\Pi_{li}^\nu\hat{\nu}-\Pi_{li}^\eta\hat{\eta}, & i=l+1,\cdots,N
\end{array}\right.\\
&\hat{e}_{zj} =z_j-\Pi_{z j}^\nu\hat{\nu}-\Pi_{zj}^\eta \hat{\eta}, \quad j=1,\cdots,M,\\
& \hat{e}_{\bar{\eta} i}= \bar{\eta}_i-\Pi_{\bar{\eta} i}^\nu\hat{\nu}-\Pi_{\bar{\eta} i}^\eta \hat{\eta}, \quad
i=1,\cdots,l.
\end{align}
Their dynamics are governed by
\begin{align}
  \dot{\hat{e}}_{x i} &=\hat{A}_ie_{xi}-\hat{D}_i\sum_{j=1}^M h_{ij}G_j\hat{e}_{zj}+\epsilon \Pi_{fi}^{\bar{\eta}}G_B\sum_{j=1}^M h_{ij}G_j\hat{e}_{zj}, \nonumber
  \\ & \quad \quad \quad \quad \quad \quad \quad \quad \quad \quad \quad \quad  i=1,\cdots,l,\\
  \dot{\hat{e}}_{x i} &=\hat{A}_ie_{xi}-\hat{D}_i\sum_{j=1}^M h_{ij}G_j\hat{e}_{zj}, \nonumber \\ 
  &\quad \quad \quad \quad \quad \quad \quad \quad \quad \quad \quad \quad  i=l+1,\cdots,N,
  \\ \dot{\hat{e}}_{z j} &= E_j\hat{e}_{z j}+ F_j\sum_{i=1}^Nh_{ij}\hat{C}_ie_{xi}+ F_j\sum_{i=1}^l h_{ij} G_Q\hat{e}_{\bar{\eta} i}, \nonumber
  \\ &\quad \quad \quad \quad \quad \quad \quad \quad \quad \quad \quad \quad j=1,\cdots,M,\\
  \dot{\hat{e}}_{\bar{\eta}i} & = G_S\hat{e}_{\bar{\eta}i}-\epsilon G_B\sum_{j=1}^Mh_{ij}G_j\hat{e}_{zj}, \nonumber
  \\ & \quad \quad \quad \quad \quad \quad \quad \quad \quad \quad \quad \quad i=1,\cdots,l,
  \end{align}
which can be further written in the following compact form
\begin{multline}
\begin{bmatrix} \dot{\hat{e}}_x \\  \dot{\hat{e}}_z \\ \dot{\hat{e}}_\eta \end{bmatrix}
= \\
\begin{bmatrix}
\overline{AN}       & -\overline{HDG}+\epsilon\begin{bmatrix} \overline{\Pi_LB} \\ 0 \end{bmatrix}   & 0 \\
\overline{H^TFC}    & \overline{EM}                                 & \overline{H^T_LF}(I_l\otimes G_Q) \\
0                   & -\epsilon(I_l\otimes G_B)\overline{H_LG}      & I_l\otimes G_s
\end{bmatrix}
\begin{bmatrix} \hat{e}_x \\ \hat{e}_z \\ \hat{e}_\eta \end{bmatrix},
\end{multline}
where $\hat{e}_x=[\hat{e}_{x1}^T,\cdots,\hat{e}_{xN}^T]^T$, $\hat{e}_z = [\hat{e}^T_{z1}, \cdots,\hat{e}_{zM}^T]^T$,
$\hat{e}_\eta=[\hat{e}_{\eta1}^T, \cdots,\hat{e}_{\eta l}^T]^T$,
\begin{align*}
&\overline{\Pi_LB} = \begin{bmatrix}
h_{11} \Pi_{f1}^{\bar{\eta}}G_BG_1      &  \cdots       & h_{1M} \Pi_{f1}^{\bar{\eta}}G_BG_M \\
\vdots                                  &  \vdots       &  \vdots\\
h_{l1} \Pi_{fl}^{\bar{\eta}}G_BG_1      &  \cdots       & h_{lM} \Pi_{fl}^{\bar{\eta}}G_BG_M
\end{bmatrix}, 
\end{align*}
\begin{align*}
&\overline{H_LG} = \begin{bmatrix}
  h_{11}G_1         & \cdots        & h_{1M}G_M \\
  \vdots            & \vdots        & \vdots\\
  h_{l1}G_1         & \cdots        & h_{lM}G_M
\end{bmatrix}, \\
&\overline{H^T_LF} = \begin{bmatrix}
  h_{11}F_1         & \cdots        & h_{lM}F_1 \\
  \vdots            & \vdots        & \vdots\\
  h_{1M}F_M         & \cdots        & h_{lM}F_M
\end{bmatrix}.
\end{align*}
Using the technique similar to that used for the proof of Theorem \ref{th02} and using Lemma \ref{le01}, it can be derived that there is a positive scalar $\epsilon^*$ such that for all $0<\epsilon<\epsilon^*$ the above system is exponentially stable, that is, all error vectors $\hat{e}_{xi}$, $\hat{e}_{zj}$ and $\hat{e}_{\bar{\eta} i}$ will converge to zero. According to \eqref{A39} and \eqref{A44},  it follows that $ v_i\rightarrow Q_v\nu_i$ for the slave nodes and $ y_j\rightarrow Q\eta_j$ for the master nodes.
\end{IEEEproof}

\bibliographystyle{unsrt}
\bibliography{/Users/yanjun/JabRef/XJ2013}
\end{document}